\documentclass[10pt,a4paper]{amsart}
\usepackage{amsfonts,layout}

\newcount\minute
\newcount\hour
\def\currenttime{%
	\minute\time
	\hour\minute
	\divide\hour60
	\the\hour:\multiply\hour60\advance\minute-\hour\the\minute}
\def\draftnote{{\it \today \quad  \currenttime \hfill  tex-file :   \jobname}}

\catcode`@=11
\@addtoreset{equation}{section}

\catcode`@=12
\newtheorem{Theorem}{Theorem}[section]

\newtheorem{Lemma}{Lemma}[section]

\newtheorem{Remark}{Remark}[section]

\newtheorem{Hyp.}{Hyp.}[section]

\begin{document}

\title[]{Precise estimates for biorthogonal families under asymptotic gap conditions}

\author{P. Cannarsa} 
\address{Dipartimento di Matematica, Universit\`a di Roma "Tor Vergata",
Via della Ricerca Scientifica, 00133 Roma, Italy}
\email{cannarsa@mat.uniroma2.it}

\author{P. Martinez} 
\address{Institut de Math\'ematiques de Toulouse, UMR CNRS 5219, Universit\'e Paul Sabatier Toulouse III \\ 118 route de Narbonne, 31 062 Toulouse Cedex 4, France} \email{Patrick.Martinez@math.univ-toulouse.fr}

\author{J. Vancostenoble} 
\address{Institut de Math\'ematiques de Toulouse, UMR CNRS 5219, Universit\'e Paul Sabatier Toulouse III \\ 118 route de Narbonne, 31 062 Toulouse Cedex 4, France} \email{Judith.Vancostenoble@math.univ-toulouse.fr}

\thanks{This research was partly supported by the Institut Mathematique de Toulouse and Istituto Nazionale di Alta Matematica through funds provided by the national group GNAMPA and the GDRE CONEDP}

\begin{abstract}
A classical and useful way to study controllability problems is the moment method developed by Fattorini-Russell \cite{FR1, FR2}, and based on the construction of suitable biorthogonal families. Several recent problems exhibit the same behaviour: the eigenvalues of the problem satisfy a uniform but rather 'bad' gap condition, and a rather 'good' but only asymptotic one. 
The goal of this work is to obtain general and precise upper and lower bounds
for biorthogonal families under these two gap conditions, and so to measure the influence of the 'bad' gap condition and the good influence of the 'good' asymptotic one. To achieve our goals, we extend some of the general results
of Fattorini-Russell \cite{FR1, FR2} concerning biorthogonal families, using complex analysis techniques developed by Seidman \cite{Seid-Avdon}, G\"uichal \cite{Guichal}, Tenenbaum-Tucsnak \cite{Tucsnak} and Lissy \cite{Lissy1, Lissy2}.  

\end{abstract}
\maketitle

\section{Introduction}
\label{sec-intro}

\subsection{Presentation of the subject} \hfill

Biorthogonal families are a classical tool in analysis.  In particular, they play a crucial role  in the so-called moment method, which was developed by Fattorini-Russell~\cite{FR1, FR2} to study controllability for parabolic equations.

Given any sequence of nonnegative real numbers, $(\lambda_n)_{n\geq 1}$, we recall that a sequence $(\sigma _m)_{m\geq 1}$ is biorthogonal to the sequence 
$(e^{\lambda _n t})_{n\geq 1}$ in $L^2(0,T)$ if

$$ \forall m, n \geq 1, \quad \int _0 ^T \sigma _m (t) e^{\lambda _n t } \, dt = \begin{cases} 1 \text{ if } m=n \\ 0 \text{ if } m\neq n \end{cases} .$$
The goal of this paper is to provide explicit and precise upper and lower bounds for the biorthogonal family $(\sigma _m)_{m\geq 1}$ under the following gap conditions:
\begin{itemize}
\item a `global gap condition':
\begin{equation}
\label{hyp-gap}
\forall n \geq 1, \quad 0 < \gamma _{min} \leq \sqrt{\lambda _{n+1}} - \sqrt{\lambda _n} \leq \gamma _{max} ,
\end{equation}
\item and an `asymptotic gap condition':
\begin{equation}
\label{hyp-gap-asympt}
\forall n \geq N^*, \quad \gamma _{min} ^* \leq \sqrt{\lambda _{n+1}} - \sqrt{\lambda _n} \leq \gamma _{max}^*,
\end{equation}
where $\gamma _{max}^*-\gamma _{min}^*<\gamma _{max}-\gamma _{min}$.
\end{itemize}
Before explaining why we are interested in such a question, let us describe some of the main results of the literature on this subject.


\subsection{The context} \hfill

Among the most important applications of biorthogonal families to control theory are those to the  null controllability and sensitivity of control costs to  parameters. Major contributions in such directions are the following:
\begin{itemize}
\item Fattorini-Russell \cite{FR1, FR2}, Hansen \cite{Hansen},  and Ammar Khodja-Benabdallah-Gonz\'alez Burgos-de Teresa \cite{Assia1} studied the existence of  biorthogonal sequences and their application to controllability for various equations;
\item for nondegenerate parabolic equations and dispersive equations, Seidman \cite{Seidman84}, G\"uichal \cite{Guichal}, Seidman-Avdonin-Ivanov \cite{Seid-Avdon}, Miller \cite{Miller}, Tenen\-baum-Tucsnak \cite{Tucsnak}, and Lissy \cite{Lissy1, Lissy2} studied the dependence of the null controllability cost $C_T$   with respect to the time $T$ (as $T \to 0$, the so-called 'fast control problem') and with respect to the domain, obtaining extremely sharp estimates of the constants $c(\Omega)$ and $C(\Omega)$ that appear in
$$ e^{c(\Omega)/T} \leq C_T \leq e^{C(\Omega)/T};$$
\item Coron-Guerrero \cite{Coron}, Glass \cite{Glass}, Lissy \cite{Lissy2} investigated the vanishing viscosity problem:
$$ \begin{cases}
y_t + M y_x - \varepsilon y_{xx} = 0 , \quad x\in (0,L), \\
y(0,t) = f(t) ,
\end{cases} $$ 
obtaining sharp estimates of the null controllability cost with respect to the time $T$,  the transport coefficient $M$,  the size of the domain $L$, and  the diffusion coefficient $\varepsilon$;
\item  in \cite{memoire, CMV-cost-weak}, we studied  the dependence of the controllability cost with respect to the degeneracy parameter $\alpha$ for the  degenerate parabolic equation
$$ u_t - (x^\alpha u_x)_x = 0 , \quad x\in (0,\ell) .$$
\end{itemize}

There is a common feature in these works: they depend on some parameter $p$, and this parameter forces the eigenvalues to satisfy \eqref{hyp-gap} (sometimes after normalization)  with gap bounds
$\gamma _{min} (p)$ and $\gamma _{max}(p)$ such that
$$ \gamma _{min}(p) \to 0 \quad \text{ and/or} \quad \gamma _{max} (p)\to \infty .$$
This fact makes it necessary to have general and precise estimates with respect to the main parameters that appear in the problem.

In \cite{CMV-cost-weak}, we proved the following general result: given $T>0$ and a family $(\lambda _n)_{n\geq 1}$ of nonnegative real numbers that satisfy the 'global gap condition' \eqref{hyp-gap}, then: 
\begin{itemize}
\item every family $(\sigma _m)_{m\geq 1}$, biorthogonal to $(e^{\lambda _n t})_{n\geq 1}$ in $L^2(0,T)$, satisfies the  lower estimate
\begin{equation}
\label{estim-MCRF-low}
\Vert \sigma _m \Vert _{L^2(0,T)} ^2 \geq b_m  \, e^{-2\lambda _{m} T} \, e^{\frac{1}{2T (\gamma _{\text{max}}) ^2}} ,
\end{equation}
with an explicit value of $b_m = b_m (T,\gamma _{max},m)$ (rational in $T$);
\item  there exists a family $(\sigma _m)_{m\geq 1}$, biorthogonal to $(e^{\lambda _n t})_{n\geq 1}$ in $L^2(0,T)$, that satisfies the  upper estimate
\begin{equation}
\label{estim-MCRF-up}
\Vert \sigma _m \Vert _{L^2(0,T)} ^2 
\leq B_m  \, e^{-2\lambda _{m} T} e^{\frac{C_u}{T (\gamma _{min}  )^2}} e^{C_u \frac{\sqrt{\lambda _m}}{\gamma _{min} }} ,
\end{equation}
with an explicit value of $B_m= B_m (T,\gamma _{min},m)$ (rational in $T$).
\end{itemize}
The bounds  \eqref{estim-MCRF-low} and \eqref{estim-MCRF-up} above describe quite precisely the behavior of the biorthogonal family, in particular in short time. Estimate \eqref{estim-MCRF-up} is in the spirit of \cite{FR1,FR2} but the dependence with respect to $T$ when $T\to 0^+$ is completely explicit, and  assumption \eqref{hyp-gap} is a little more general than the asymptotic development of the eigenvalues used in Tenenbaum-Tucsnak \cite{Tucsnak} or Lissy \cite{Lissy1, Lissy2}:
$$ \lambda _n = rn^2 + O(n).$$
Moreover, \eqref{hyp-gap}  explains the role of  $\gamma _{min}$ and $\gamma _{max}$  in the analysis of the biorthogonal family:  
 $\gamma_{min}$ determines, essentially, the growth rate of the upper bound for $(\sigma _m)_{m\geq 1}$ while
  $\gamma _{max}$ gives the lower bound.


\subsection{Motivations and main results of  this paper} \hfill
 
 Even though the aforementioned results give a fairly good picture of the properties of the family $(\sigma _m)_m$, 
 some delicate issues remain to be analysed and will be addressed in this paper. For instance, one would like to understand the dependence of 
 the family $(\sigma _m)_m$ with respect to relevant parameters that come into play. Typical examples of such problems are the following ones.
 \begin{itemize}
\item For the 1D degenerate parabolic equation
$$ u_t - (x^\alpha u_x)_x = 0 , \quad x\in (0,\ell),$$
 the eigenvalues $\lambda _{\alpha,n}$ of the associated elliptic operator (with suitable boundary conditions) can be expressed using the zeros of Bessel functions (\cite{Gueye}) and depend on the degeneracy parameter $\alpha \in (0,2)$.  One can then prove (see \cite{memoire, CMV-cost-loc}) that  the global gap condition \eqref {hyp-gap} is satisfied only with
$$ \gamma _{max} (\alpha) \geq c (2-\alpha ) ^{2/3} , $$
with $c>0$, while the asymptotic gap condition \eqref{hyp-gap-asympt} is satisfied with
$$ \gamma _{max} ^* (\alpha) \leq c^* (2-\alpha ),$$
where $c^* >0$, after the rank
$$ N^* (\alpha) = \frac{1}{2-\alpha};$$
in this case
$$ \frac{\gamma _{max} (\alpha)}{\gamma _{max}^* (\alpha)} \longrightarrow +\infty \quad\text{ as }  \alpha \to 2 ^- ;$$ 
hence it is natural to think that the better asymptotic gap \eqref{hyp-gap-asympt} could be used to improve the estimate \eqref{estim-MCRF-low} of the associated biorthogonal sequences, but the fact that 
$$ N^* (\alpha)\longrightarrow +\infty \quad\text{ as }  \alpha \to 2 ^-$$ 
is certainly to be taken into account.

\item In  2D problems such as the Grushin equation (see \cite{Karine-Grushin1, Karine-Grushin2}),  where  the solution is decomposed into Fourier modes, one has to give uniform bounds for a certain sequence of elliptic problems,  the eigenvalues of which satisfy \eqref{hyp-gap} and \eqref{hyp-gap-asympt} with some  $\gamma _{min}(m)$, $\gamma _{min}^*(m)$ and $N^* (m)$ such that
$$ \frac{\gamma _{min} (m)}{\gamma _{min}^* (m)} \longrightarrow 0 \text{ as } m\to \infty $$
and 
$$ N^* (m)\longrightarrow +\infty \quad\text{ as } m\to \infty ;$$
once again, it is natural to think that the better asymptotic gap \eqref{hyp-gap-asympt} could be used to improve the estimate \eqref{estim-MCRF-up} of the associated biorthogonal sequence, but the fact that 
$N^* (m)\to +\infty$ as $m\to \infty$ 
is certainly to be taken into account.
\end{itemize}

The above discussion motivates the general question whether estimates \eqref{estim-MCRF-low} and \eqref{estim-MCRF-up} can be improved when \eqref{hyp-gap} is combined with the asymptotic condition \eqref{hyp-gap-asympt}. This is exactly what we prove in this paper:  roughly speaking, 
\eqref{estim-MCRF-low} and \eqref{estim-MCRF-up} hold true  replacing 
$\gamma _{min}$ by $\gamma _{min} ^*$ and $\gamma _{max}$ by $\gamma _{max}^*$.
Moreover, the fact the 'good' gap condition \eqref{hyp-gap-asympt}
holds true only after the $N^*$ first eigenvalues has a cost, and we obtain a precise estimate for that cost.
Our main results (Theorem~\ref{thm-biortho1-gen} and \ref{thm-guichal-gen}) are the following: under \eqref{hyp-gap-asympt}, we prove that:
\begin{itemize}
\item every family $(\sigma _m)_{m\geq 1}$, biorthogonal to $(e^{\lambda _n t})_{n\geq 1}$ in $L^2(0,T)$, satisfies the  lower estimate
\begin{equation}
\label{estim-MCRF-low*}
\Vert \sigma _m \Vert _{L^2(0,T)} ^2 \geq b^* _m  \, e^{-2\lambda _{m} T} \, e^{\frac{2}{T (\gamma _{max}^*) ^2}} ,
\end{equation}
where the 'cost' $b_m ^* = b^* _m (T,\gamma _{max}, \gamma _{max}^*, N^*,m)$ is a rational function of $T$ that we determine explicitly, and
\item  there exists a biorthogonal family that satisfies
\begin{equation}
\label{estim-MCRF-up*}
 \Vert \sigma _m \Vert _{L^2(0,T)} ^2 
\leq B_m ^* \, e^{-2\lambda _{m} T} e^{\frac{C_0}{T (\gamma _{min} ^* )^2}} e^{C_0 \frac{\sqrt{\lambda _m}}{\gamma _{min} ^*}} ,
\end{equation}
where  $C_0 >0$ is a universal constant and $B^* _m (T,\gamma _{min}, \gamma _{min}^*, N^*,m)$ is  a rational function of $T$ that we determine explicitly.
\end{itemize}
Let us observe that the presence of the exponential factors $e^{\frac{2}{T (\gamma _{max}^*) ^2}}$ and $e^{\frac{C_0}{T (\gamma _{min} ^* )^2}}$ in  \eqref{estim-MCRF-low*} and \eqref{estim-MCRF-up*}
is quite natural and has already been pointed out by Seidman-Avdonin-Ivanov~\cite{Seid-Avdon}, Tenenbaum-Tucsnak~\cite{Tucsnak}, and Lissy~\cite{Lissy1, Lissy2}
(see also Haraux~\cite{Haraux} and Komornik~\cite{Kom-Lor} for a closely related context).  On the other hand, the precise estimate of the behavior of $b_m ^*$ and $B_m ^*$ with respect to  parameters, that we develop in this paper, is completely new and will be crucial for the sensitivity analysis of  control costs to be performed in \cite{CMV-cost-loc}.

Our proofs are based on complex analysis techniques and Hilbert space methods developed by Seidman-Avdonin-Ivanov~\cite{Seid-Avdon} and G\"uichal~\cite{Guichal}. We have also used an idea from Tenenbaum-Tucsnak~\cite{Tucsnak} and Lissy~\cite{Lissy1, Lissy2}, based on the introduction of an extra parameter depending on $T$ and the gap conditions. 



\subsection{Plan of the paper} \hfill

The paper is organized as follows: 
\begin{itemize}
\item in section \ref{s-results}, we state  our results;
\item section \ref{sec-biortho} is devoted to the proof of Theorem \ref{thm-biortho1-gen} (construction of a biorthogonal family and derivation of upper bounds);
\item section \ref{sec-guichal} is devoted to the proof of Theorem \ref{thm-guichal-gen} (lower bounds for biorthogonal families).
\end{itemize}


\section{Setting of the problem and main results}
\label{s-results}

\subsection{Existence of a suitable biorthogonal family and upper bounds} \hfill
\label{s5}

We will establish the following results, that in some sense provide a more precise version of properties observed by Fattorini and Russell \cite{FR1, FR2} (in short time),  much in the spirit of Tenenbaum-Tucsnak \cite{Tucsnak} and Lissy \cite{Lissy1, Lissy2} (with a slightly weakened assumption on the eigenvalues).

\begin{Theorem}
\label{thm-biortho1-gen}
Assume that 
$$ \forall n\geq 1, \quad \lambda_n \geq 0, $$
and that there is some $0 < \gamma _{\text{min}} < \gamma _{\text{min}} ^*$ such that
\begin{equation}
\label{gap-min}
\forall n \geq 1, \quad \sqrt{\lambda _{n+1}} - \sqrt{\lambda _{n}}  \geq \gamma _{\text{min}} ,
\end{equation}
and
\begin{equation}
\label{gap-min*}
\forall n \geq N^*, \quad \sqrt{\lambda _{n+1}} - \sqrt{\lambda _{n}}  \geq \gamma _{\text{min}}^* .
\end{equation}
Denote
\begin{equation}
\label{def-M*}
M^* := (1-\frac{\gamma _{min}}{\gamma _{min} ^*}) (N^* -1) .
\end{equation}
Then there exists a family $(\sigma _{m} ^+)_{m\geq 1}$ which is biorthogonal to the family $(e^{\lambda _{n}t})_{n\geq 1}$ in $L^2(0,T)$:
\begin{equation}
\label{*famillebi_qm-gen}
\forall m,n \geq 1, \quad \int _0 ^T \sigma _{m} ^+ (t)  e^{\lambda _{n}t} \, dt = \delta _{mn} .
\end{equation}
Moreover, it satisfies: there is some universal constant $C$ independent of $T$, $\gamma _{\text{min}}$, $\gamma _{\text{min}}^*$, $N^*$ and $m$ such that, for all $m\geq 1$, we have
\begin{equation}
\label{*famillebi_qm-norme-gen}
\Vert \sigma _{m} ^+ \Vert _{L^2(0,T)} ^2  
\leq e^{-2\lambda _{m} T} e^{\frac{C}{T (\gamma _{min} ^* )^2}}
e^{C \frac{\sqrt{\lambda _m}}{\gamma _{min} ^*}} 
 B^* (T,\gamma _{min},\gamma _{min}^*, N^*,m )  ,
\end{equation}
where


\begin{multline}
\label{eq(B>}
B^* (T,\gamma _{min},\gamma _{min}^*, N^*,m )
\\
= \begin{cases}
C_u \Bigl(
\frac{(8M^*)!}{(\lambda _{m}  (\gamma _{min} ^*)^2 T^2)^{4M^*}} +1 \Bigr)
e ^{C_u M^*}
e^{C_u \frac{\lambda _{N^*}}{\gamma _{min} \sqrt{\lambda _m}}}
(\frac{1}{T^{3/2}} + \frac{1}{(\gamma _{min} ^* )^2 T^2})
\text{ if } T \leq \frac{1}{(\gamma _{min} ^*) ^2}
\\
C_u \Bigl(
(\frac{(\gamma _{min} ^*)^{8M^*} (8M^*)!}{\lambda _{m} ^{4M^*}}) +1\Bigr)
e ^{C_u M^*}
e^{C_u \frac{\lambda _{N^*}}{\gamma _{min} \sqrt{\lambda _m}}}
((\gamma _{min} ^*) ^2 + (\gamma _{min} ^*) ^3) 
\text{ if } T \geq \frac{1}{(\gamma _{min} ^*) ^2}
\end{cases}  .
\end{multline}


\end{Theorem}

\begin{Remark} 
{\rm Theorem \ref{thm-biortho1-gen} completes and improves several earlier results, in particular Theorem 1.5 of Fattorini-Russell \cite{FR2} and \cite{CMV-cost-weak}, providing the explicit dependence of the $L^2$ bound with respect to $\gamma _{min}$, $\gamma _{min}^*$  in short time.
It is useful in several problems, in which $\gamma _{min} \to 0$ with respect to some parameter, which occurs is several cases, see, e.g. \cite{FR3}, \cite{Karine-Grushin1}.
We will apply the construction used by Seidman, Avdonin and Ivanov in \cite{Seid-Avdon}, which has the advantage to be completely explicit (which is not the case for the construction of \cite{FR1,FR2, FR3, Hansen, Assia1}, since there is a contradiction argument), combined with some ideas coming from the construction of Tenenbaum-Tucsnak \cite{Tucsnak} and Lissy \cite{Lissy1}, adding some parameter, in order to obtain precise results.}
\end{Remark}


\subsection{General lower bounds} \hfill

We generalise a result by G\"uichal \cite{Guichal} to prove the following

\begin{Theorem}
\label{thm-guichal-gen}
Assume that 
$$ \forall n\geq 1, \quad \lambda_n \geq 0, $$
and that there are $0<\gamma _{min} \leq \gamma _{\text{max}}^* \leq \gamma _{\text{max}}$ such that
\begin{equation}
\label{gap-max}
\forall n \geq 1, \quad \gamma _{\text{min}} \leq \sqrt{\lambda _{n+1}} - \sqrt{\lambda _{n}}  \leq \gamma _{\text{max}} ,
\end{equation}
and
\begin{equation}
\label{gap-max*}
\forall n \geq N_*, \quad \sqrt{\lambda _{n+1}} - \sqrt{\lambda _{n}}  \leq \gamma _{\text{max}}^* .
\end{equation}
Then any family $(\sigma _{m} ^+)_{m\geq 1}$ which is biorthogonal to the family $(e^{\lambda _{n}t})_{n\geq 1}$ in $L^2(0,T)$ (hence that satisfies \eqref{*famillebi_qm-gen})
satisfies:
\begin{equation}
\label{csq-gap-max}
\Vert \sigma _{m} ^+ \Vert _{L^2(0,T)} ^2
\geq e^{-2\lambda _{m} T} \, e^{\frac{2}{T (\gamma _{\text{max}}^*) ^2}} \, b^* (T,\gamma_{max},\gamma_{max} ^*, N_*, \lambda _1,m)^2 ,
\end{equation}
where $b^*$ is rational in $T$ (and explictly given in the key Lemma \ref{lem-dist-gapmax}).
\end{Theorem}

\begin{Remark} 
{\rm Theorem \ref{thm-guichal-gen} completes a result of G\"uichal \cite{Guichal} and is useful in several problems, in which $\gamma _{max} \to \infty$ with respect to some parameter, which occurs is several cases, see, e.g. \cite{FR3}, and \cite{CMV-cost-loc}.
It is to be noted that the behaviour with respect to $m$ can perhaps be improved,
comparing with Theorem 1.1 of Hansen \cite{Hansen}. It would be interesting to investigate this.
}
\end{Remark}


\section{Proof of Theorem \ref{thm-biortho1-gen}} 
\label{sec-biortho}


\subsection{The general strategy} \hfill

It begins with the following remarks: if the family $(\sigma _{m} ^+)_{m\geq 1}$ is biorthogonal to the family $(e^{\lambda _{n}t})_{n\geq 1}$, then 
$$ \forall m,n \geq 1, \quad  \int _0 ^T \sigma _{m} ^+ (T-t) e^{\lambda _{n}(T-t)} \, dt = \delta _{mn} ,$$
hence
$$ \forall m,n \geq 1, \quad  \int _0 ^T \Bigl( \sigma _{m} ^+ (T-t) e^{\lambda _{m}T } \Bigr) e^{-\lambda _{n}t } \, dt = \delta _{mn} ,$$
hence the family $(s_{m})_{m\geq 0}$ defined by
$$ s_{m}(t) := \sigma _{m} ^+ (T-t) e^{\lambda _{m}T } $$
is biorthogonal to the family $(e^{-\lambda _{n}t})_{n\geq 1}$ in $L^2(0,T)$. Now extend $s_{m}$ by $0$ outside $(0,T)$, and consider its Fourier transform
$$ \forall z\in \Bbb C, \quad \mathcal{F}(s_{m}) (z) := \int _{\Bbb R} s_{m} (t) e^{-izt} \, dt .$$
For all $m\geq 1$, $\mathcal{F}(s_{m})$ is the Fourier transform of a compactly supported function, hence it is an entire function over $\Bbb C$, and it satisfies
$$ \forall m,n \geq 1, \quad  \mathcal{F}(s_{m}) (-i \lambda _{n})= \delta _{mn} ,$$
and it is of exponential type:
$$ \vert \mathcal{F}(s_{m}) (z) \vert \leq \Bigl( \int _0 ^T \vert s_{m} (t) \vert \, dt \Bigr) e^{T \vert z \vert }; $$
and also
$$  \mathcal{F}(s_{m}) (z) = \int _{-T/2} ^{T/2} s_{m} (\tau+ \frac{T}{2}) e^{-iz(\tau + \frac{T}{2})} \, d \tau ,$$
hence
$$  \mathcal{F}(s_{m}) (-z) e^{-iz\frac{T}{2}}  = \int _{-T/2} ^{T/2} s_{m} (\tau+ \frac{T}{2}) e^{iz\tau} \, d\tau ,$$
and
$$ \vert \mathcal{F}(s_{m}) (-z) e^{-iz \frac{T}{2}} \vert \leq \Bigl( \int _{-T/2} ^{T/2} \vert s_{m} (\tau+ \frac{T}{2}) \vert \, d\tau \Bigr) e^{\frac{T}{2}\vert z \vert} .$$

Now we recall the Paley-Wiener theorem (\cite{Young}): if $f: \Bbb C \to \Bbb C$ is an entire function 
of exponential type, such that there exist nonnegative constants $C,A$ such that
$$ \forall z\in \Bbb C, \quad \vert f(z) \vert \leq C e^{A\vert z \vert },$$
and if $f \in L^2 (\Bbb R)$, then there exists $\phi \in L^2(-A,A)$ such that 
$$ f(z) = \int _{-A} ^A \phi (\tau) e^{iz\tau} \, d\tau .$$
One of the objects of \cite{Seid-Avdon} is to prove the existence of a sequence $(f_m)_m$ of entire functions satisfying
\begin{equation}
\label{*cond-gene}
\begin{cases}
\forall m,n \geq 1, \quad f_m (-i\lambda _n) = \delta _{mn} , \\
\forall z \in \Bbb C, \quad \vert f_m (-z) e^{-iz \frac{T}{2}} \vert \leq C_m e^{\frac{T}{2}\vert z \vert} ,\\
\forall m \geq 1, \quad f_m \in L^2 (\Bbb R) 
\end{cases}
\end{equation}
(see Theorem 2 and Lemma 3 in \cite{Seid-Avdon}) under some general assumptions on the sequence $(\lambda_n)_n$. If we can apply such a result in our context (hence with our sequence $(\lambda_{n})_n$), then the two last properties together with the Paley-Wiener theorem will imply that there exists some $\phi_m \in L^2 (-\frac{T}{2}, \frac{T}{2})$ such that
$$ f_m (-z) e^{-iz \frac{T}{2}} = \int _{-T/2}^{T/2} \phi_m (\tau) e^{iz\tau} \,d\tau ,$$
hence
$$ f_m(z) = \int _0^T \phi _m (t-\frac{T}{2}) e^{-izt} \, dt ,$$
and then
$$ \int _0^T \phi _m (t-\frac{T}{2}) e^{-\lambda _{n} t} \, dt = f_m  (-i\lambda _{n}) = \delta _{mn} ,$$
hence $(\phi _m (t-\frac{T}{2}))_m$ will be biorthogonal to the family $(e^{-\lambda _{n} t})_n$, and $(\sigma _{m} ^+ (t))_m$ defined by
$$\sigma _{m} ^+ (t) = \phi _m (\frac{T}{2}-t) e^{-\lambda _{m}T } $$
will be biorthogonal to the family $(e^{\lambda _{n} t})_n$ in $L^2(0,T)$, as desired. Moreover
$$ \Vert \sigma _{m} ^+ \Vert _{L^2(0,T)} ^2 = e^{-2\lambda _{m}T } 
\int _{-T/2} ^{T/2} \phi _m (\tau) ^2 \, d\tau 
\leq C e^{-2\lambda _{m}T } \Vert f_m  \Vert _{L^2(\Bbb R)} ^2 $$
using the Parseval theorem.

Now, it remains to construct such entire functions $f_m$. The idea is to consider the natural infinite product that satisfies the first condition of \eqref{*cond-gene}, $f_m (-i\lambda _n) = \delta _{mn}$, and to multiply it by a so-called 'mollifier', in such a way that the other two conditions of \eqref{*cond-gene} will be also satisfied. Hence one has to estimate the growth of the natural infinite product, and then to choose a choose a suitable mollifier. This is what is performed in \cite{Seid-Avdon}.
For our problem, our task will be to add the dependency into the parameters $\gamma _{min}$, $\gamma _{min}^*$ and $T$, and to understand specifically the 
behaviour of the natural infinite product, the mollifier and at the end of $\Vert \sigma _{m} ^+ \Vert _{L^2(0,T)}$ with respect to $\gamma_{min}$ and $T$. 
We will modify a little the construction of \cite{Seid-Avdon}, in order to obtain optimal results in our context, see Lemma \ref{*lem-mollifier}, and specifically the definition \eqref{*defPq} of the mollifier, where the additional parameter $N'$ will be chosen of the size $\frac{1}{T(\gamma _{min}^*) ^2}$, see \eqref{*sam7oct1}.


\subsection{The counting function} \hfill

Consider 
$$ \forall \rho >0, \quad N_n (\rho) := \text{ card } \{k, 0< \vert \lambda _n - \lambda _k \vert \leq \rho \}
.$$
We prove the following:

\begin{Lemma}
\label{lem-counting}
a) Assume that the gap assumption \eqref{gap-min} is satisfied; then 
\begin{equation}
\label{counting}
\forall n\geq 0, \forall \rho >0, \quad N_n (\rho) \leq 2 \frac{\sqrt{\rho}}{\gamma _{min}}.
\end{equation}

b) Assume that the gap assumptions \eqref{gap-min}-\eqref{gap-min*} are satisfied; then

\begin{itemize}
\item when $n=N^*$: 
\begin{equation}
\label{counting*=}
\forall \rho >0, \quad N_{N^*} (\rho) \leq 
\begin{cases} 
\frac{\sqrt{\rho}}{\gamma _{min}} + \frac{\sqrt{\rho}}{\gamma _{min}^*} \quad & \text{ if } \rho \leq \lambda _{N^*}  \\
N^* -1 + \frac{\sqrt{\rho}}{\gamma _{min}^*} \quad & \text{ if } \rho \geq \lambda _{N^*} 
\end{cases} ,
\end{equation}

\item when $n>N^*$:
\begin{equation}
\label{counting*>}
\forall n > N^*, \forall \rho >0, \quad 
N_n (\rho) \leq 
\begin{cases} 
\frac{2\sqrt{\rho}}{\gamma _{min}^*} \quad & \text{ if } \rho \leq \lambda _{n} -\lambda _{N^*}  
\\
\frac{\sqrt{\rho}}{\gamma _{min}} + \frac{\sqrt{\rho}}{\gamma _{min}^*} \quad & \text{ if } \lambda _{n} -\lambda _{N^*} \leq \rho \leq \lambda _{n}  
\\
n -1 +\frac{\sqrt{\rho}}{\gamma _{min}^*} \quad & \text{ if } \rho \geq \lambda _n  
\end{cases} ,
\end{equation}

\item when $n<N^*$
\begin{itemize}
\item when $\lambda _n \leq \lambda _{N^*}-\lambda _{n}$, then
\begin{equation}
\label{counting*<}
\forall \rho >0, \quad 
N_n (\rho) \leq 
\begin{cases} 
\frac{2\sqrt{\rho}}{\gamma _{min}} \quad & \text{ if } \rho \leq \lambda _{n} 
\\
n-1 + \frac{\sqrt{\rho}}{\gamma _{min}} \quad & \text{ if } \lambda _{n} \leq \rho \leq \lambda _{N^*}-\lambda _{n} \\
N^* -1 + \frac{\sqrt{\rho}}{\gamma _{min}^*} \quad & \text{ if } \rho \geq \lambda _{N^*} -\lambda _{n}
\end{cases} , 
\end{equation}
\item when $\lambda _n \geq \lambda _{N^*}-\lambda _{n}$, then
\begin{equation}
\label{counting*<2}
\forall \rho >0, \quad 
N_n (\rho) \leq 
\begin{cases} 
\frac{2\sqrt{\rho}}{\gamma _{min}} \quad & \text{ if } \rho \leq \lambda _{N^*}-\lambda _{n} \\
N^*-n + \frac{\sqrt{\rho}}{\gamma _{min}} + \frac{\sqrt{\rho}}{\gamma _{min}^*} \quad & \text{ if } \lambda _{N^*}-\lambda _{n} \leq \rho \leq \lambda _{n} \\
N^* -1 + \frac{\sqrt{\rho}}{\gamma _{min}^*} \quad & \text{ if } \rho \geq \lambda _{n}
\end{cases} . 
\end{equation}
\end{itemize}
\end{itemize}

\end{Lemma}

\noindent \begin{Remark} The main point in \eqref{counting*=}-\eqref{counting*<2} is to observe that $N_n$ behaves as $\frac{\sqrt{\rho}}{\gamma _{min}^*}$, as $\rho \to +\infty$, and to compute all the needed additional constants. 
\end{Remark}

\noindent {\it Proof of Lemma \ref{lem-counting}.}
Take $k>n$. Then 
$$ \lambda _k - \lambda _n = \sqrt{\lambda _k} ^2 - \sqrt{\lambda _n}  ^2 = (\sqrt{\lambda _k}  - \sqrt{\lambda _n} )(\sqrt{\lambda _k}  + \sqrt{\lambda _n} ),$$
and the gap assumption \eqref{gap-min} insures that
$$ \sqrt{\lambda _k}  - \sqrt{\lambda _n}  \geq (k-n) \gamma _{min} , \quad \sqrt{\lambda _k} + \sqrt{\lambda _n}  \geq (k-n) \gamma _{min} + 2 \sqrt{\lambda _n}  ,$$
hence
$$ \lambda _k - \lambda _n  \geq (k-n)^2 \gamma _{min}^2 ,$$
and
$$ k>n \text{ and } \lambda _k - \lambda _n \leq \rho \quad \implies  \quad k-n \leq \frac{\sqrt{\rho}}{\gamma _{min}} .$$
Similarly, 
$$ k< n \text{ and } \lambda _n - \lambda _k \leq \rho \quad \implies  \quad n-k \leq \frac{\sqrt{\rho}}{\gamma _{min}} .$$
Hence
$$ N_n (\rho) \leq 2 \frac{\sqrt{\rho}}{\gamma _{min}}. $$
This proves \eqref{counting}.

Now we prove \eqref{counting*=}-\eqref{counting*<2}: let us introduce 
$$ N_n ^+ (\rho) := \text{ card } \{k >n, \lambda _k - \lambda _n \leq \rho \},
\quad 
N_n ^- (\rho) := \text{ card } \{k <n, \lambda _n - \lambda _k \leq \rho \} .$$
We distinguish the three cases.

\begin{itemize}
\item When $n=N^*$: from the previous study, we see that
$$ N_{N^*} ^+ (\rho) \leq \frac{\sqrt{\rho}}{\gamma _{min}^*} ,
\quad \text{ and } \quad N_{N^*} ^- (\rho) \leq 
\begin{cases} 
\frac{\sqrt{\rho}}{\gamma _{min}} \quad & \text{ if } \rho \leq \lambda _{N^*}  \\
N^* -1 \quad & \text{ if } \rho \geq \lambda _{N^*} 
\end{cases} ;
$$
this gives that
$$ N_{N^*} (\rho) \leq 
\begin{cases} 
\frac{\sqrt{\rho}}{\gamma _{min}} + \frac{\sqrt{\rho}}{\gamma _{min}^*} \quad  & \text{ if } \rho \leq \lambda _{N^*}  \\
N^* -1 + \frac{\sqrt{\rho}}{\gamma _{min}^*} \quad & \text{ if } \rho \geq \lambda _{N^*} 
\end{cases} ,
$$
which gives \eqref{counting*=}.

\item When $n>N^*$: now we have
$$ N_n ^+ (\rho) \leq \frac{\sqrt{\rho}}{\gamma _{min}^*} ,$$
and
$$  N_n ^- (\rho) \leq 
\begin{cases} 
\frac{\sqrt{\rho}}{\gamma _{min}^*} \quad & \text{ if } \rho \leq \lambda _{n} -\lambda _{N^*}  
\\
\frac{\sqrt{\rho}}{\gamma _{min}}  \quad & \text{ if } \lambda _{n} -\lambda _{N^*} \leq \rho \leq \lambda _{n}  \\
n -1 \quad & \text{ if } \rho \geq \lambda _n 
\end{cases} ,
$$
which gives \eqref{counting*>}. 

\item When $n<N^*$: now we have
$$ N_n ^- (\rho) \leq 
\begin{cases} 
\frac{\sqrt{\rho}}{\gamma _{min}} \quad & \text{ if } \rho \leq \lambda _{n}  \\
n -1 \quad & \text{ if } \rho \geq \lambda _n 
\end{cases} ,
$$
and
$$ N_n ^+ (\rho) \leq 
\begin{cases} 
\frac{\sqrt{\rho}}{\gamma _{min}}  \quad  & \text{ if } \rho \leq \lambda _{N^*}-\lambda _{n}  \\
N^*-n + \frac{\sqrt{\rho}}{\gamma _{min}^*} \quad & \text{ if } \rho \geq \lambda _{N^*}-\lambda _n 
\end{cases} ,
$$
hence when $\lambda _n \leq \lambda _{N^*}-\lambda _{n}$ we have
$$ N_n (\rho) \leq 
\begin{cases} 
\frac{2\sqrt{\rho}}{\gamma _{min}} \quad & \text{ if } \rho \leq \lambda _{n} \\
n-1+ \frac{\sqrt{\rho}}{\gamma _{min}} \quad & \text{ if } \lambda _{n} \leq \rho \leq \lambda _{N^*}-\lambda _{n} \\
N^*-1 + \frac{\sqrt{\rho}}{\gamma _{min}^*} \quad & \text{ if } \rho \geq \lambda _{N^*} -\lambda _{n}
\end{cases} , $$
which gives \eqref{counting*<}, and similar estimates when $\lambda _n \geq \lambda _{N^*}-\lambda _{n}$, which give \eqref{counting*<2}. \qed

\end{itemize}

Before going further, let us give another estimate of the counting function,
which reveals to be more practical and more natural, since it gives a better understanding of the role
of the different parameters:

\begin{Lemma}
\label{lem-counting-new}
Assume that the gap assumptions \eqref{gap-min}-\eqref{gap-min*} are satisfied; then

\begin{itemize}
\item when $n=N^*$: 
\begin{equation}
\label{counting*=-new}
\forall \rho >0, \quad N_{N^*} (\rho) \leq 
\begin{cases} 
\frac{\sqrt{\rho}}{\gamma _{min}} + \frac{\sqrt{\rho}}{\gamma _{min}^*} \quad & \text{ if } \rho \leq \lambda _{N^*}  \\
(1-\frac{\gamma _{min}}{\gamma _{min} ^*}) (N^*-1) + 2 \frac{\sqrt{\rho}}{\gamma _{min}^*} \quad & \text{ if } \rho \geq \lambda _{N^*} 
\end{cases} ,
\end{equation}

\item when $n>N^*$:
\begin{equation}
\label{counting*>-new}
\forall \rho >0, \quad 
N_n (\rho) \leq 
\begin{cases} 
\frac{2\sqrt{\rho}}{\gamma _{min}^*} \quad & \text{ if } \rho \leq \lambda _{n} -\lambda _{N^*}  
\\
\frac{\sqrt{\rho}}{\gamma _{min}} + \frac{\sqrt{\rho}}{\gamma _{min}^*} \quad & \text{ if } \lambda _{n} -\lambda _{N^*} \leq \rho \leq \lambda _{n}  \\
(1-\frac{\gamma _{min}}{\gamma _{min} ^*})(N^* -1) +2\frac{\sqrt{\rho}}{\gamma _{min}^*} \quad & \text{ if } \rho \geq \lambda _n  
\end{cases} ,
\end{equation}

\item when $n<N^*$:
\begin{equation}
\label{counting*<-new}
\forall \rho >0, \quad N_n (\rho) \leq 
\begin{cases} 
\frac{2\sqrt{\rho}}{\gamma _{min}} \quad & \text{ if } \rho \leq \max \{\lambda _{n} , \lambda _{N^*}  - \lambda _n \}
\\
(1-\frac{\gamma _{min}}{\gamma _{min} ^*})(N^* -1) + (1+\sqrt{2}) \frac{\sqrt{\rho}}{\gamma _{min}^*} \quad & \text{ if } \rho \geq \max \{\lambda _{n} , \lambda _{N^*}  - \lambda _n \} 
\end{cases}  ,
\end{equation}
and also
\begin{equation}
\label{counting*<-new2}
\forall \rho >0, \quad N_n (\rho) \leq 
\begin{cases} 
\frac{2\sqrt{\rho}}{\gamma _{min}} \quad & \text{ if } \rho \leq \lambda _{N^*} \\
(1-\frac{\gamma _{min}}{\gamma _{min} ^*})(N^*-1) + 2 \frac{\sqrt{\rho}}{\gamma _{min}^*} \quad & \text{ if } \rho \geq \lambda _{N^*} 
\end{cases} .
\end{equation}

\end{itemize}

\begin{Remark} {\rm Lemma \ref{lem-counting-new} enlightens the role of the quantity $(1-\frac{\gamma _{min}}{\gamma _{min} ^*})(N^*-1)$ (denoted $M^*$ in \eqref{def-M*}); when 
$\gamma _{min} =\gamma _{min} ^*$ or if $N^*=1$, this quantity is equal to zero, and we logically find estimates similar to the ones of Lemma \ref{lem-counting} (i.e. the "1 gap condition"); in the more interesting case where 
$\gamma _{min} < \gamma _{min} ^*$ and $N^*>1$, this quantity measures the increase of the counting function with respect to the "1 gap condition".

Let us note also that we expect that \eqref{counting*<-new} holds true with $2$ instead of $1+\sqrt{2}$, however we could not prove it in full generality.
}
\end{Remark}

\end{Lemma}

\noindent {\it Proof of Lemma \ref{lem-counting-new}.}
\begin{itemize}
\item When $n=N^*$, it is sufficient to note that
$$ \sqrt{\lambda _{N^*}} \geq \gamma _{min} (N^* -1),$$
hence, when $\rho \geq \lambda _{N^*}$, we have
\begin{multline*}
N^*-1 = (1-\frac{\gamma _{min}}{\gamma _{min} ^*}) (N^*-1) + \frac{\gamma _{min}}{\gamma _{min} ^*} (N^*-1)
\\
\leq (1-\frac{\gamma _{min}}{\gamma _{min} ^*}) (N^*-1) + \frac{\sqrt{\lambda _{N^*}}}{\gamma _{min} ^*}
\leq (1-\frac{\gamma _{min}}{\gamma _{min} ^*}) (N^*-1) + \frac{\sqrt{\rho}}{\gamma _{min} ^*} ;
\end{multline*}
this and \eqref{counting*=} imply \eqref{counting*=-new}.

\item When $n>N^*$: when $\rho \geq \lambda _n$, we have
$$ \sqrt{\lambda _n} \geq \sqrt{\lambda _{N^*}} + (n-N^*) \gamma _{min} ^*
\geq (N^*-1) \gamma _{min} + (n-N^*) \gamma _{min} ^* ,$$
hence
$$ n -1 \leq \frac{\sqrt{\lambda _n}}{\gamma _{min} ^*} + (1-\frac{\gamma _{min}}{\gamma _{min} ^*}) (N^*-1)
\leq \frac{\sqrt{\rho}}{\gamma _{min} ^*} + (1-\frac{\gamma _{min}}{\gamma _{min} ^*}) (N^*-1) ;$$
this estimate and \eqref{counting*>} imply \eqref{counting*>-new}.

\item When $n<N^*$, we obtain \eqref{counting*<-new2} proceeding in the same way: 
when $\rho \leq \lambda _{N^*}$, then clearly $N_n (\rho)$ is less than the number of terms that would be at both sides, for which the gap of their square root would be $\gamma _{min}$, hence 
$$ N_n (\rho) \leq 2\frac{\sqrt{\rho}}{\gamma _{min}} ;$$
when $\rho \geq \lambda _{N^*}$, then clearly one has all the $N^*-1$ first terms, and the others, for which the gap of their square root is $\gamma _{min}^*$, hence
$$ N_n (\rho) \leq N^*-1 + \frac{\sqrt{\rho}}{\gamma _{min}^*} ;$$
but then
\begin{multline*}
 N^*-1 = (1- \frac{\gamma _{min}}{\gamma _{min}^*}) (N^*-1) + \frac{\gamma _{min}}{\gamma _{min}^*} (N^*-1)
\\
\leq (1- \frac{\gamma _{min}}{\gamma _{min}^*}) (N^*-1) + \frac{\sqrt{\lambda _{N^*}}}{\gamma _{min}^*} \leq (1- \frac{\gamma _{min}}{\gamma _{min}^*}) (N^*-1) + \frac{\sqrt{\rho}}{\gamma _{min}^*} ,
\end{multline*}
which gives \eqref{counting*<-new2};

\item finally we prove \eqref{counting*<-new}: in the same way, if $\rho \leq \max \{\lambda _{n} , \lambda _{N^*}  - \lambda _n \}$ one has immediately
 $$ N_n (\rho) \leq \frac{2\sqrt{\rho}}{\gamma _{min}} ;$$
 when $\rho \geq \max \{\lambda _{n} , \lambda _{N^*}  - \lambda _n \}$, then 
we already know from \eqref{counting*<} and \eqref{counting*<2} that
$$ N_n (\rho) \leq N^* -1 + \frac{\sqrt{\rho}}{\gamma _{min} ^*},$$
hence
$$ N_n (\rho) \leq (1-\frac{\gamma _{min}}{\gamma _{min} ^*}) (N^*-1) + \frac{\gamma _{min}}{\gamma _{min} ^*} (N^*-1)+ \frac{\sqrt{\rho}}{\gamma _{min} ^*};$$
since
$$ \sqrt{\lambda _{N^*}} \geq \gamma _{min} (N^*-1) $$
and
$$ \rho \geq \max \{\lambda _{n} , \lambda _{N^*}  - \lambda _n \} \quad \implies \quad
\rho \geq \frac{1}{2} \Bigl( \lambda _{n}+ (\lambda _{N^*}  - \lambda _n) \Bigr) 
= \frac{1}{2} \lambda _{N^*} ,$$
we deduce that
$$ \gamma _{min} (N^*-1) \leq \sqrt{\lambda _{N^*}} \leq \sqrt{2\rho} ,$$
hence
$$  N_n (\rho) \leq (1-\frac{\gamma _{min}}{\gamma _{min} ^*}) (N^*-1) + (1+\sqrt{2}) \frac{\sqrt{\rho}}{\gamma _{min} ^*} ,$$
which is \eqref{counting*<-new}.

\end{itemize}

\noindent This concludes the proof of Lemma \ref{lem-counting-new}. \qed


\subsection{A Weierstrass product} \hfill

 Motivated by \cite{Seid-Avdon}, we consider
\begin{equation}
\label{*-prodinf}
F_m (z) := 
\prod _{k=1, k\neq m} ^\infty \Bigl( 1 - \Bigl( \frac{iz - \lambda _{m}}{\lambda _{k} - \lambda _{m}}\Bigr) ^2 \Bigr) .
\end{equation}
Then the growth in $k$ of $\lambda _{k}$ ensures that this infinite product converges uniformly over all the compact sets, hence $F_m$ is well-defined and entire over $\Bbb C$. Moreover
$$ F_m (-i\lambda _{n}) = \prod _{k=1, k\neq m} ^\infty \Bigl( 1 - \Bigl( \frac{\lambda _{n} - \lambda _{m}}{\lambda _{k} - \lambda _{m}}\Bigr) ^2 \Bigr) 
= \begin{cases} 0 \quad \text{ if } m\neq n , \\ 1  \quad \text{ if } m = n ,\end{cases},$$
hence 
\begin{equation}
\label{*-prodinf-cond1}
\forall m, n \geq 1, \quad F_m (-i\lambda _{n}) = \delta _{mn} .
\end{equation}
We are going to estimate the growth of $F_m$. We prove the following

\begin{Lemma}
\label{*lem-growthFm}
a) Assume that the gap assumption \eqref{gap-min} is satisfied. Then the function $F_m$ satisfies the following growth estimate: there is some uniform constant $C_u$ (independent of $m$, $\gamma _{min}$, and $z$) such that
\begin{equation}
\label{*eq-growthFm}
\forall z\in \Bbb C, \quad \vert F_m (z) \vert 
\leq e^{\frac{C_u}{\gamma _{min}} \sqrt{\lambda _{m} }}
 e^{\frac{C_u}{\gamma _{min}} \sqrt{\vert z\vert} } .
\end{equation} 

b) Assume that the gap assumptions \eqref{gap-min}-\eqref{gap-min*} are satisfied. Then the function $F_m$ satisfies the following growth estimate: there is some uniform constant $C_u$ (independent of $m$, $\gamma _{min}$, $\gamma _{min} ^*$ and $z$), such that
\begin{equation}
\label{*eq-growthFm-gene}
\forall m \geq 1, \forall z \in \Bbb C, \quad \vert F_m (z) \vert 
\leq B_m \, q_m (\vert z \vert) \, e^{C_u \frac{\sqrt{\vert z \vert}}{\gamma _{min}^* }} ,
\end{equation} 
with
\begin{equation}
\label{*eq-growthFm-Bm-unif}
\forall m\geq 1, \quad 
B_m \leq 
e^{ \frac{C_u}{\gamma _{min} } \frac{\lambda _{N^*}}{\sqrt{\lambda _m }}}
\, e^{\frac{C_u}{\gamma _{min}^*} \sqrt{\lambda _m }}
\end{equation}
and
\begin{equation}
\label{*eq-growthFm-qm-unif}
\forall m\geq 1, \quad
q_m (\vert x \vert ) \leq 
\Bigl( 3+2\frac{\vert z \vert ^2}{\lambda _m ^2} \Bigr) ^{M^*} 
,
\end{equation}
where $M^*$ has been defined in \eqref{def-M*}.
\end{Lemma}

\noindent \begin{Remark}
{\rm The main point in Lemma \ref{*lem-growthFm} is to obtain estimates of the growth of $F_m$ in
$e^{C_u \frac{\sqrt{\vert z \vert}}{\gamma _{min}^* }}$ under \eqref{gap-min}-\eqref{gap-min*},
with explicit constants (given in \eqref{*eq-growthFm-Bm-unif} and \eqref{*eq-growthFm-qm-unif}), that will help us in the following.
Comparing with \eqref{*eq-growthFm}, this gives a better idea of the improvement brought by 'large' gap $\gamma_{min}^*$ and of the price to pay due to the 'small' gap $\gamma_{min}$ for the $N^*$ first eigenvalues. In fact we will first prove the following better estimates: \eqref{*eq-growthFm-gene} holds true with

\begin{equation}
\label{*eq-growthFm-Bm}
B_m = 
\begin{cases}
e^{(\frac{8}{\gamma _{min}}+\frac{C_u}{\gamma _{min}^*})\sqrt{\lambda _m}+ \frac{8}{\gamma _{min}}\sqrt{\lambda _{N^*}-\lambda _m} } \quad & \text{ if } m<N^* \\
e^{(\frac{4}{\gamma _{min}}+ \frac{4+C_u}{\gamma _{min}^*}) \sqrt{\lambda _{N^*} }} \quad & \text{ if } m=N^* \\
e^{(\frac{4}{\gamma _{min}}+ \frac{8+C_u}{\gamma _{min}^*}) \sqrt{\lambda _{m} } - \frac{4}{\gamma _{min}} \sqrt{\lambda _{m} - \lambda _{N^*}} }   \quad & \text{ if } m>N^* 
\end{cases}
\end{equation} 
and
\begin{equation}
\label{*eq-growthFm-qm}
q_m (\vert z \vert)=
\begin{cases}
\Bigl( 1+2 \frac{\vert z \vert ^2 + \lambda _m ^2}{\max (\lambda _m  , \lambda _{N^*}-\lambda _m) ^2} \Bigr)^{M^*} \quad &\text{ if } m<N^* 
\\
\Bigl( 1+ 2\frac{\vert z \vert ^2 + (\lambda _{N^*}) ^2}{(\lambda _{N^*}) ^2} \Bigr) ^{M^*}  \quad & \text{ if } m=N^* 
\\ 
\Bigl( 1+2\frac{\vert z \vert ^2 + \lambda _m ^2}{\lambda _m ^2} \Bigr) ^{M^*}
\quad & \text{ if } m>N^* 
\end{cases} ,
\end{equation}
and this easily implies \eqref{*eq-growthFm-Bm-unif} and \eqref{*eq-growthFm-qm-unif}.
}
\end{Remark}

{\it Proof of Lemma \ref{*lem-growthFm}.} Note that
$$
F_m (z-i \lambda _{m}) 
= \prod _{k=1, k\neq m} ^\infty \Bigl( 1 + \frac{z^2 }{(\lambda _{k} - \lambda _{m})^2} \Bigr) , 
$$
hence (following \cite{Seid-Avdon})

\begin{multline*}
\ln \vert F_m (z-i \lambda _{m}) \vert 
= \sum _{k=1, k\neq m} ^\infty \ln \Bigl \vert 1 + \frac{z^2 }{(\lambda _{k} - \lambda _{m})^2} \Bigr\vert 
\leq \sum _{k=1, k\neq m} ^\infty \ln \Bigl( 1 + \frac{\vert z \vert ^2 }{(\lambda _{k} - \lambda _{m})^2} \Bigr) 
\\
= \int _0 ^\infty \ln \Bigl( 1 + \frac{\vert z \vert ^2 }{\rho ^2} \Bigr)  d N_{m} (\rho)
= 2 \int _0 ^\infty \frac{N_{m} (\rho)}{\rho} \frac{\vert z \vert ^2 }{\vert z \vert ^2 + \rho ^2}  \, d\rho 
\\
\end{multline*}
Then we distinguish several cases:

\begin{itemize}
\item Under only \eqref{gap-min} we deduce from \eqref{counting} that

$$ 2 \int _0 ^\infty \frac{N_{m} (\rho)}{\rho} \frac{\vert z \vert ^2 }{\vert z \vert ^2 + \rho ^2}  \, d\rho 
\leq \frac{4}{\gamma _{min}}  \int _0 ^\infty \frac{1}{\sqrt{\rho}} \frac{\vert z \vert ^2 }{\vert z \vert ^2 + \rho ^2}  \, d\rho 
=  \Bigl( \frac{4}{\gamma _{min}} \int _0 ^\infty \frac{1}{\sqrt{ s}} \frac{1}{1 + s^2}  \, ds \Bigr) \sqrt{\vert z \vert } .
$$
Then changing $z-i \lambda _{m}$ into $z$, 
$$
\ln \vert F_m (z) \vert 
\leq \Bigl( \frac{4}{\gamma _{min}} \int _0 ^\infty \frac{1}{\sqrt{ s}} \frac{1}{1 + s^2}  \, ds \Bigr)  ( \sqrt{\vert z\vert}  + \sqrt{\lambda _{m} } ), 
$$
which gives \eqref{*eq-growthFm}.

\item Under \eqref{gap-min}-\eqref{gap-min*} and when $m = N^*$, we derive from \eqref{counting*=-new} that
\begin{multline*}
2 \int _0 ^\infty \frac{N_{N^*} (\rho)}{\rho} \frac{\vert z \vert ^2 }{\vert z \vert ^2 + \rho ^2}  \, d\rho 
\\
= 2 \int _0 ^{\lambda _{N^*}} \frac{N_{N^*} (\rho)}{\rho} \frac{\vert z \vert ^2 }{\vert z \vert ^2 + \rho ^2}  \, d\rho + 2 \int _{\lambda _{N^*}} ^\infty \frac{N_{N^*} (\rho)}{\rho} \frac{\vert z \vert ^2 }{\vert z \vert ^2 + \rho ^2}  \, d\rho
\\
\leq 2 \int _0 ^{\lambda _{N^*}} \sqrt{\rho} (\frac{1}{\gamma _{min}} + \frac{1}{\gamma _{min}^*}) \frac{\vert z \vert ^2 }{\rho (\vert z \vert ^2 + \rho ^2)}  \, d\rho
+ 2 \int _{\lambda _{N^*}} ^\infty (M^* + \frac{2\sqrt{\rho}}{\gamma _{min}^*}) \frac{\vert z \vert ^2 }{\rho(\vert z \vert ^2 + \rho ^2)}  \, d\rho 
\\
\leq 2 (\frac{1}{\gamma _{min}}+ \frac{1}{\gamma _{min}^*}) \sqrt{\vert z \vert} \int _0 ^{\lambda _{N^*}/\vert z \vert}  \frac{1 }{ \sqrt{s} (1 +  s ^2)}  \, ds
\\
+ 2 M^* \int _{\lambda _{N^*}/\vert z \vert} ^\infty \frac{1}{s (1 +  s ^2)}  \,  ds
+ \frac{4\sqrt{\vert z \vert}}{\gamma _{min}^*} \int _{\lambda _{N^*}/\vert z \vert} ^\infty   \frac{1 }{\sqrt{s} (1 +  s ^2)}  \,  ds
\\
\leq 2(\frac{1}{\gamma _{min}}+ \frac{1}{\gamma _{min}^*}) \sqrt{\vert z \vert}  \frac{2 \sqrt{\lambda _{N^*}}}{\sqrt{\vert z \vert}} + 2 M^* [\ln \frac{s}{\sqrt{1+s^2}} ]_{\lambda _{N^*}/\vert z \vert} ^\infty
+ \frac{4\sqrt{\vert z \vert}}{\gamma _{min}^*} \int _{0} ^\infty   \frac{1 }{\sqrt{s} (1 +  s ^2)}  \,  ds
\\
\leq 4 \sqrt{\lambda _{N^*}}(\frac{1}{\gamma _{min}}+ \frac{1}{\gamma _{min}^*}) + M^* \ln (1+ \frac{\vert z \vert ^2}{(\lambda _{N^*}) ^2})+ c_u \frac{\sqrt{\vert z \vert}}{\gamma _{min}^*} .
\end{multline*}
Then changing $z-i \lambda _{N^*}$ into $z$, 
$$
\ln \vert F_{N^*} (z) \vert 
\leq 4 \sqrt{\lambda _{N^*}}(\frac{1}{\gamma _{min}}+ \frac{1}{\gamma _{min}^*}) + M^* \ln (1+ 2\frac{\vert z \vert ^2 + \lambda _{N^*} ^2}{(\lambda _{N^*}) ^2})+ c_u \frac{\sqrt{\vert z \vert}+ \sqrt{\lambda _{N^*} }}{\gamma _{min}^*}, 
$$
which gives \eqref{*eq-growthFm-gene} with the $B_m$ and $q_m$ given in \eqref{*eq-growthFm-Bm} and \eqref{*eq-growthFm-qm}.

\item Under \eqref{gap-min}-\eqref{gap-min*} and when $m>N^*$, applying the same method, we derive from \eqref{counting*>-new} that
\begin{multline*}
2 \int _0 ^\infty \frac{N_{m} (\rho)}{\rho} \frac{\vert z \vert ^2 }{\vert z \vert ^2 + \rho ^2}  \, d\rho
\leq 4(\frac{1}{\gamma _{min}}+\frac{2}{\gamma _{min}^*}) \sqrt{\lambda _m}
- \frac{4}{\gamma _{min}}\sqrt{\lambda _m- \lambda _{N^*}} 
\\
+ M^* \ln (1+\frac{\vert z \vert ^2}{\lambda _m ^2})
+ c_u \frac{\sqrt{\vert z \vert}}{\gamma _{min}^*} .
\end{multline*}
Then changing $z-i \lambda _{m}$ into $z$, we obtain \eqref{*eq-growthFm-gene} with the related $B_m$ and $q_m$  given in \eqref{*eq-growthFm-Bm} and \eqref{*eq-growthFm-qm}. 

\item Under \eqref{gap-min}-\eqref{gap-min*} and when $m<N^*$, applying the same method, we derive from \eqref{counting*<-new} that
\begin{multline*}
2 \int _0 ^\infty \frac{N_{m} (\rho)}{\rho} \frac{\vert z \vert ^2 }{\vert z \vert ^2 + \rho ^2}  \, d\rho
\leq \frac{8}{\gamma _{min}}(\sqrt{\lambda _m}+ \sqrt{\lambda _{N^*} - \lambda _m}) 
\\
+ M^* \ln (1 +  \frac{\vert z \vert ^2}{\max \{ \lambda _m , \lambda _{N^*}-\lambda _m \} ^2})
+ c_u \frac{\sqrt{\vert z \vert}}{\gamma _{min}^*} .
\end{multline*}
Then changing $z-i \lambda _{m}$ into $z$, we obtain \eqref{*eq-growthFm-gene} with the related $B_m$ and $q_m$  given in \eqref{*eq-growthFm-Bm} and \eqref{*eq-growthFm-qm}.  \qed

\end{itemize}


\subsection{A suitable mollifier} \hfill

Motivated by \cite{Seid-Avdon}, we made in \cite{CMV-cost-weak} the following construction:
consider $T'>0$, $N' \geq 1$, $a_k := \frac{C_{N',T'}}{k^2}$ with
$$ C_{N',T'} := \frac{T'}{2 \sum _{k=N'} ^\infty \frac{1}{k^2} },$$
in order that
$$ \sum _{k=N'} ^\infty a_k = \frac{T'}{2},$$
and finally
\begin{equation}
\label{*defPq}
P_{N',T'} (z) := e^{iz\frac{T'}{2}} \prod _{k=N'} ^\infty \cos (a_k z) .
\end{equation}
Then we have the following

\begin{Lemma}(\cite{CMV-cost-weak})
\label{*lem-mollifier}

\begin{enumerate}
\item The regularity and the growth of $P_{N',T'}$ over $\Bbb C$: The function $P_{N',T'}$ is entire over $\Bbb C$ and satisfies
\begin{equation} 
\label{*-mol-cond1}
\begin{cases}
P_{N',T'}(0)=1, \\
\forall z\in \Bbb C \text{ such that } \Im z \geq 0, \quad \vert P_{N',T'}(z) \vert \leq 1,\\
\forall z\in \Bbb C , \quad \vert e^{-iz\frac{T'}{2}} P_{N',T'}(z) \vert \leq e^{\vert z \vert \frac{T'}{2}} .
\end{cases} 
\end{equation}

\item The behaviour of $P_{N',T'}$ over $\Bbb R$: there exist $\theta _0 >0$, $\theta _1 >0$, both independent of $N'$ and $T'$ such that $P_{N',T'}$ satisfies
\begin{equation} 
\label{*-mol-cond2}
\begin{cases}
\Bigl( \frac{C_{N',T'} \vert x \vert }{\theta _0} \Bigr) ^{1/2} + 1 \geq N' \implies
 \ln \vert P_{N',T'}(x) \vert \leq -\frac{\theta _1}{2^3} \Bigl( \frac{C_{N',T'} \vert x \vert}{\theta _0} \Bigr)^{1/2}  , \\
\Bigl( \frac{C_{N',T'} \vert x \vert}{\theta _0} \Bigr) ^{1/2} + 1 \leq N' \implies
\ln \vert P_{N',T'}(x) \vert \leq -\frac{\theta _1}{(N')^3}  \Bigl( \frac{C_{N',T'} \vert x \vert}{\theta _0} \Bigr)^2 .
\end{cases}
\end{equation}

\item The behaviour of $P_{N',T'}$ over $i\Bbb R _+$: there is some constant $\theta _2 >0$, independent of $N'$ and $T'$, such that $P_{N',T'}$ satisfies
\begin{equation} 
\label{*-mol-cond3}
\forall x \in \Bbb R _+, \quad P_{N',T'}(ix) \geq  e^{-\theta _2 \sqrt{C_{N',T'}x } } .
\end{equation}
 \end{enumerate}

\end{Lemma}

\noindent The Proof of Lemma \ref{*lem-mollifier} follows by elementary analysis techniques.
In the following we are going to use the mollifier $P_{N',T'}$ to construct the biorthogonal family. 


\subsection{A sequence of holomorphic functions satisfying \eqref{*cond-gene}} \hfill

Consider
\begin{equation}
\label{*-deff_m}
\forall m \geq 0, \forall z \in \Bbb C, \quad 
f_{m,N',T'} (z) :=  F_m (z) \frac{P_{N',T'} (-z)}{P_{N',T'} (i\lambda _{m})} .
\end{equation}
We will make the following choices:
\begin{itemize}
\item for $T'$:
\begin{equation}
\label{eq-T'}
T' := \min \{ T, \frac{1}{(\gamma _{min} ^*)^2} \};
\end{equation}
\item for $N'$: we choose it such that
\begin{equation}
\label{*sam7oct1-new}
N' \geq 2+ \frac{\theta_ 3}{(\gamma _{min}^*)^2 T'}
\end{equation}
with a suitable $\theta _3$ (independent of $T>0$ and of $m\geq 0$, and given in \eqref{*ven6oct1}).
\end{itemize}
Then we will prove the following

\begin{Lemma}
\label{lem-propfmN}
When $T'$ and $N'$ satisfy \eqref{eq-T'} and \eqref{*sam7oct1-new}, the functions $f_{m,N',T'}$ are entire and satisfy the following properties:
\begin{itemize}
\item for all $m,n \geq 1$, we have
\begin{equation}
\label{*-deff_m-cond1}
f_{m,N',T'}(-i \lambda _n) = \delta _{mn} ;
\end{equation}
\item for all $m \geq 1$, for all $\varepsilon >0$, there exists $C_{m,N',T',\varepsilon}>0$ such that
\begin{equation}
\label{*-deff_m-cond2}
\forall z \in \Bbb C, \quad 
\vert f_{m,N',T'} (-z)  e^{-iz\frac{T}{2}} \vert 
\leq C_{m,N',T',\varepsilon} e^{(\frac{T}{2} + \varepsilon) \vert z \vert } ;
\end{equation}
\item  for all $m \geq 1$, $f_{m,N',T'} \in L^2 (\Bbb R)$.
\end{itemize}

\end{Lemma}
Then we will be in position to apply the Paley-Wiener theorem and to construct the desired biorthogonal sequence.

\noindent {\it Proof of Lemma \ref{lem-propfmN}.}
First, the function $f_{m,N',T'}$ is well-defined since $P_{N',T'} >0$ on $i\Bbb R_+$, and is entire since $F_m$ and $P_{N',T'}$ are entire.
Next, using \eqref{*-prodinf-cond1}, we have \eqref{*-deff_m-cond1}.
Next, concerning the exponential type: using \eqref{*eq-growthFm-gene} and \eqref{*-mol-cond1}, we have
\begin{multline*}
 \vert f_{m,N',T'} (-z)  e^{-iz\frac{T}{2}} \vert 
= \vert F_m (-z)  \vert \, \vert P_{N',T'} (z)  e^{-iz\frac{T'}{2}} \vert \, \vert e^{-iz\frac{T-T'}{2}} \vert \, \frac{1}{P_{N',T'} (i\lambda _{m})},$$
\\
\leq \frac{1}{P_{N',T'} (i\lambda _{m})}B_m q_m(\vert z \vert) e^{C_u\frac{\sqrt{\vert z\vert}}{\gamma _{min}^*} } e^{\vert z \vert \frac{T'}{2}} e^{\vert z \vert \frac{T-T'}{2}};
\end{multline*}
but for all $\varepsilon >0$ we have
$$ C_u\frac{\sqrt{\vert z\vert}}{\gamma _{min}^* } 
= C_u\frac{\sqrt{\varepsilon\vert z\vert}}{\gamma _{min}^*\sqrt{\varepsilon }} 
\leq \frac{C_u ^2}{2 (\gamma _{min} ^*)^2\varepsilon } + \frac{\varepsilon }{2} \vert z\vert , $$
and 
$$ q_m(\vert z \vert) \leq c_m ' e^{\frac{\varepsilon }{2} \vert z\vert} $$
which imply \eqref{*-deff_m-cond2}.
Finally, concerning the behaviour over $\Bbb R$, we deduce from \eqref{*eq-growthFm}, \eqref{*-mol-cond2} and \eqref{*-mol-cond3} that, if $\vert x \vert$ is large enough, then
$$ \vert f_{m,N',T'} (x) \vert \leq  \frac{1}{P_{N',T'} (i\lambda _{m})} B_m q_m (\vert x \vert) e^{\frac{C_u}{\gamma _{min}^*} \sqrt{\vert x\vert} }
e^{-\frac{\theta_1}{8} \Bigl( \frac{C_{N',T'} \vert x \vert}{\theta _0} \Bigr)^{1/2}} ,$$
hence $ f_{m,N',T'} \in L^2 (\Bbb R)$ if 
$$ \frac{C_u}{\gamma _{min}^*} -\frac{\theta_1}{8} \Bigl( \frac{C_{N',T'} }{\theta _0} \Bigr)^{1/2} < 0 ,$$
which is true choosing $T'$ and $N'$ satisfying \eqref{eq-T'} and \eqref{*sam7oct1-new}: indeed, 
$$ C_{N',T'} = \frac{T'}{2 \sum _{k=N'} ^\infty  \frac{1}{k^2} } ,$$
and
$$  \frac{1}{N'} = \int _{N'} ^\infty \frac{1}{y^2} \, dy 
\leq  \sum _{k=N'} ^\infty \frac{1}{k^2} \leq \int _{N'-1} ^\infty \frac{1}{y^2} \, dy = \frac{1}{N'-1} ,$$
hence
\begin{equation}
\label{estim-CNT}
\frac{(N'-1)T'}{2} \leq C_{N',T'} \leq \frac{N'T'}{2} .
\end{equation}
Hence, if
\begin{equation}
\label{*ven6oct1}
(N'-1)T' > \frac{\theta _3}{(\gamma _{min} ^*)^2} \quad \text{ with } \quad \theta _3 := \frac{2^7\theta _0 C_u ^2}{\theta _1 ^2},
\end{equation}
we obtain that $ f_{m,N',T'} \in L^2 (\Bbb R)$. And one easily verifies that $T'$, $N'$ satisfying \eqref{eq-T'} and \eqref{*sam7oct1-new} satisfy also \eqref{*ven6oct1}. This completes the proof of Lemma \ref{lem-propfmN}. \qed


\subsection{The resulting biorthogonal sequence} \hfill

With our choices, the function $x\mapsto f_{m,N',T'} (-x) e^{-ixT/2}$ is in $L^2 (\Bbb R)$, and we can consider its Fourier transform $\phi _{m,N',T'}$:
$$ \phi _{m,N',T'} (\xi) := \frac{1}{2\pi} \int _{\Bbb R} f_{m,N',T'}(-x) e^{-ix \frac{T}{2}} e^{-i\xi x } \, dx .$$
It  is well-defined since $f_{m,N',T'} \in L^2(\Bbb R)$, and the Paley-Wiener theorem (\cite{Young} p. 100) shows that 
$ \phi _{m,N',T'}$ is compactly supported in $[-\frac{T}{2} - \varepsilon, \frac{T}{2} + \varepsilon]$ (thanks to \eqref{*-deff_m-cond2}). Since this is true for all $\varepsilon >0$, $ \phi _{m,N',T'}$ is compactly supported in $[-\frac{T}{2}, \frac{T}{2}]$. 

To obtain good results, we will choose $N'$ satisfying the stronger property:
\begin{equation}
\label{*sam7oct1}
2+ \frac{\theta_ 3}{(\gamma _{min}^*)^2 T'} \leq N' \leq 4+ \frac{\theta_ 3}{(\gamma _{min}^*)^2 T'}.
\end{equation}

Then we have the following

\begin{Lemma}
\label{*lem-propf_m}
Take $T'$ and $N'$ satisfying \eqref{eq-T'} and \eqref{*sam7oct1}, and consider
\begin{equation}
\label{def-famillebiortho}
 \sigma _{m,N',T'} ^+ (t) := \phi _{m,N',T'} (\frac{T}{2}-t) e^{-\lambda _{m} T}.
 \end{equation}
Then the family $(\sigma _{m,N',T'} ^+)_{m\geq 1}$ is biorthogonal to the family $(e^{\lambda _{n}t})_{n\geq 1}$ in $L^2(0,T)$:
\begin{equation}
\label{*famillebi_qm}
\forall m,n \geq 1, \quad \int _0 ^T \sigma _{m,N',T'} ^+ (t)  e^{\lambda _{n}t} \, dt = \delta _{mn} .
\end{equation}
Moreover, it satisfies: there is some universal constant $C_u$ independent of $T$, $\gamma _{min}$, $\gamma _{min} ^*$, $N^*$ and $m$ such that, for all $m\geq 1$, we have
\begin{equation}
\label{*famillebi_qm-norme}
\Vert \sigma _{m,N',T'} ^+ \Vert _{L^2(0,T)} ^2  
\leq C_u e^{-2\lambda _{m} T} 
e^{C_u \frac{\sqrt{\lambda _m}}{\gamma _{min} ^*}} 
 B(T,\gamma _{min},\gamma _{min}^*, N^*,m ) ,
\end{equation}
where $B(T,\gamma _{min},\gamma _{min}^*, N^*,m )$ is given by \eqref{eq(B>}.

\end{Lemma}

\noindent {\it Proof of Lemma \ref{*lem-propf_m}.}  The Fourier inversion theorem gives that

$$
f_{m,N',T'}(-x) e^{-ix \frac{T}{2}} = \int _{\Bbb R}  \phi _{m,N',T'} (\xi) e^{i\xi x } \, d\xi = \int _{-T/2} ^{T/2}  \phi _{m,N',T'} (\xi) e^{i\xi x } \, d\xi .
$$
Then
\begin{multline*}
\int _0 ^T \sigma _{m,N',T'} ^+ (t)  e^{\lambda _{n}t} \, dt
= \int _0 ^T \phi _{m,N',T'} (\frac{T}{2}-t) e^{-\lambda _{m} T}  e^{\lambda _{n}t} \, dt
\\
= e^{-\lambda _{m} T} \int _{-T/2} ^{T/2} \phi _{m,N',T'} (\xi) e^{\lambda _{n} (\frac{T}{2} - \xi)} \, d\xi
= e^{-\lambda _{m} T} e^{\lambda _{n} \frac{T}{2}} \int _{-T/2} ^{T/2} \phi _{m,N',T'} (\xi) e^{-\lambda _{n} \xi} \, d\xi
\\
= e^{-\lambda _{m} T} e^{\lambda _{n} \frac{T}{2}} f_{m,N',T'}(-i\lambda _n) e^{\lambda _n \frac{T}{2}}
= f_{m,N',T'}(-i\lambda _n) e^{(\lambda _n -\lambda _{m}) T} = \delta _{mn} .
\end{multline*}
This gives \eqref{*famillebi_qm}. Concerning \eqref{*famillebi_qm-norme}, we note that the Parseval equality gives

\begin{multline}
\label{inv-Parseval}
\int _{\Bbb R} \vert f_{m,N',T'}(x) \vert ^2 \, dx 
= \int _{\Bbb R} \vert  f_{m,N',T'}(-x) e^{-ix \frac{T}{2}} \vert ^2 \, dx 
\\
= 2\pi  \int _{\Bbb R}  \vert \phi _{m,N',T'} (\xi) \vert ^2  \, d\xi 
= 2\pi  \int _{-T/2} ^{T/2} \vert \phi _{m,N',T'} (\xi) \vert ^2  \, d\xi  .
\end{multline}
Hence
$$ \Vert \sigma _{m,N',T'} ^+ \Vert _{L^2(0,T)} ^2 = e^{-2\lambda _{m} T} \int _{-T/2} ^{T/2} \vert \phi _{m,N',T'} (\xi) \vert ^2  \, d\xi = \frac{1}{2\pi} e^{-2\lambda _{m} T} \int _{\Bbb R} \vert f_{m,N',T'}(x) \vert ^2 \, dx  .$$
We need to estimate precisely the last integral. Denote
$$ X_{N',T'} := \frac{\theta _0 (N'-1)^2}{C_{N',T'}} .$$
Using \eqref{*eq-growthFm}, \eqref{*-mol-cond2} and \eqref{*-mol-cond3}, we have
\begin{multline*}
\int _{\Bbb R} \vert f_{m,N',T'}(x) \vert ^2 \, dx
= \int _{\vert x \vert \leq X_{N',T'} } \vert f_{m,N',T'}(x) \vert ^2 \, dx
+ \int _{\vert x \vert \geq X_{N',T'} } \vert f_{m,N',T'}(x) \vert ^2 \, dx
\\
\leq 2 e^{2\theta _2 \sqrt{C_{N',T'} \lambda _m}} B_m ^2
\Bigl( \int _0 ^{X_{N',T'}} q_m (x) ^2 \, e^{\frac{2C_u}{\gamma _{min}^*} \sqrt{x}} \, e^{-\frac{2\theta _1}{(N')^3} (\frac{C_{N',T'}x}{\theta _0})^2} \, dx
\\
+ \int _{X_{N',T'}} ^\infty q_m (x)^2 \,  e^{\frac{2C_u}{\gamma _{min}^*} \sqrt{x}} \, e^{-\frac{2\theta _1}{2^3} (\frac{C_{N',T'}x}{\theta _0})^{1/2}} \, dx \Bigr)
= : I^{(<)} _m+ I^{(>)}_m .
\end{multline*}
First we estimate $I^{(<)}_m$; we denote $\theta _i$ various constants independent of all the other parameters; we have
\begin{multline*}
\int _0 ^{X_{N',T'}} q_m (x) ^2 e^{\frac{2C_u}{\gamma _{min}^*} \sqrt{x}} e^{-\frac{2\theta _1}{(N')^3} (\frac{C_{N',T'}x}{\theta _0})^2} \, dx
\\
\leq q_m (X_{N',T'}) ^2 e^{\frac{2C_u}{\gamma _{min}^*} \sqrt{X_{N',T'}}} \int _0 ^{\infty} e^{-\frac{2\theta _1}{(N')^3} (\frac{C_{N',T'}x}{\theta _0})^2} \, dx
\\
\leq C'' q_m (X_{N',T'}) ^2 e^{\frac{2C_u}{\gamma _{min}^*} \sqrt{X_{N',T'}}} \frac{(N')^{3/2}}{C_{N',T'}}
\\
\leq C''' q_m (X_{N',T'}) ^2 e^{\frac{2C_u}{\gamma _{min}^*} \sqrt{X_{N',T'}}}  \Bigl( \frac{1}{T'} + \frac{1}{(T')^{3/2} \gamma _{min}^*} \Bigr) .
\end{multline*}
Using \eqref{eq-T'}, \eqref{estim-CNT} and \eqref{*sam7oct1}, we have
$$ X_{N',T'} \leq \theta _4 (\frac{1}{T'} + \frac{1}{(\gamma _{min} ^*)^2 (T')^2}),$$
$$ \frac{\sqrt{X_{N',T'}}}{\gamma _{min}^*} \leq \theta _5 (1 + \frac{1}{(\gamma _{min} ^*)^2 T'}) ,
\quad \text{ and } \quad \sqrt{C_{N',T'}} \leq  \frac{\theta _5}{\gamma _{min} ^*} ;$$
hence
$$ I^{(<)}_m \leq c_u e^{\theta _6 \frac{\sqrt{\lambda _m}}{\gamma _{min} ^*}} B_m ^2
q_m (\theta _4 (\frac{1}{T'} + \frac{1}{(\gamma _{min} ^*)^2(T')^2})) ^2 e^{\frac{\theta _5}{(\gamma _{min} ^*)^2 T'}}  \Bigl( \frac{1}{T'} + \frac{1}{(T')^{3/2} \gamma _{min}^*} \Bigr) .$$
To conclude, we will use the following basic remark: 
$$ y\in [0,1] \implies (1+y)^n\leq 2^n, \quad \text{ and } \quad y\geq 1 \implies
(1+y)^n = y^n (1+\frac{1}{y})^n \leq 2^n y^n ,$$
hence
$$ y\geq 0 \implies (1+y)^n \leq 2^n (1+ y^n) ,
 \quad \text{ and } \quad a,b \geq 0 \implies (a+b)^n \leq 2^n (a^n+b^n) .$$
Since (from \eqref{eq-T'})
$$ \frac{1}{T'} \leq \frac{1}{(\gamma _{min} ^*)^2(T')^2} ,$$
we obtain that: 
%
\begin{multline*}
q_m (\theta _4 (\frac{1}{T'} + \frac{1}{(\gamma _{min} ^*)^2(T')^2})) ^2
\leq \Bigl( 3 + \frac{2}{\lambda _m ^2} (\frac{2\theta _4 }{(\gamma _{min} ^*)^2(T')^2})^2 \Bigr) ^{2M^*}
\\
\leq 2^{2M^*} \Bigl( 3^{2M^*} + (\frac{2}{\lambda _m ^2} (\frac{2\theta _4 }{(\gamma _{min} ^*)^2(T')^2})^2) ^{2M^*} \Bigr)
\\
\leq C_u ^{2M^*} \Bigl(1 + 
(\frac{1}{\lambda _m  (\gamma _{min} ^*)^2 (T')^2})^{4M^*} \Bigr) ,
\end{multline*}
then
\begin{multline*}
\forall m\geq N^*, \quad I^{(<)}_m \leq \frac{c_u}{(T')^{3/2} \gamma _{min}^*} e^{\theta _7 \frac{\sqrt{\lambda _m}}{\gamma _{min} ^*}}
e^{\frac{\theta _5}{(\gamma _{min} ^*)^2 T'}} 
\\
e^{\theta _7 \frac{\lambda _{N^*}}{\gamma _{min} \sqrt{\lambda _m}}}
C_u ^{2M^*} \Bigl(1 + 
(\frac{1}{\lambda _m  (\gamma _{min} ^*)^2 (T')^2})^{4M^*} \Bigr) 
.
\end{multline*}

Next we estimate $I^{(>)}_m$. 
Denote
$$ L := \frac{2\theta _1}{2^3} (\frac{C_{N',T'}}{\theta _0})^{1/2}- \frac{2C_u}{\gamma _{min}^*} .
$$
One can easily check that
$$ \frac{1}{L} \leq \frac{C_u}{T' \gamma _{min} ^*} .$$
Then
\begin{multline*}
I_m ^{(>)} = 2 e^{2\theta _2 \sqrt{C_{N',T'} \lambda _m}} B_m ^2 \int _{X_{N',T'}} ^\infty q_m(x)^2 \, e^{\frac{2C_u}{\gamma _{min}^*} \sqrt{x}} \, e^{-\frac{2\theta _1}{2^3} (\frac{C_{N',T'}x}{\theta _0})^{1/2}} \, dx
\\
= 2 e^{2\theta _2 \sqrt{C_{N',T'} \lambda _m}} B_m ^2 \int _{X_{N',T'}} ^\infty q_m(x)^2 \, e^{-L \sqrt{x}} \, dx 
\\
\leq 
2 e^{\theta _6 \frac{\sqrt{\lambda _m}}{\gamma _{min} ^*}} 
e^{\theta _7 \frac{\lambda _{N^*}}{\gamma _{min} \sqrt{\lambda _m}}}
\int _0 ^\infty q_m(x)^2 \, e^{-L \sqrt{x}} \, dx 
\\
= 2 e^{\theta _6 \frac{\sqrt{\lambda _m}}{\gamma _{min} ^*}} 
e^{\theta _7 \frac{\lambda _{N^*}}{\gamma _{min} \sqrt{\lambda _m}}}
\frac{2}{L^2} \int _0 ^\infty q_m(\frac{t^2}{L^2})^2 \, e^{-t} \, dt . 
\end{multline*}
Recalling that
$$ \int _0 ^\infty t^{k} e^{-t} \, dt = k ! ,$$
we obtain
%

\begin{multline*}
I_m ^{(>)} \frac{L^2}{4} e^{-\theta _6 \frac{\sqrt{\lambda _m}}{\gamma _{min} ^*}} 
e^{-\theta _7 \frac{\lambda _{N^*}}{\gamma _{min} \sqrt{\lambda _m}}}
\leq 
2^{2M^*} \int _0 ^\infty (3^{2M^*} + \Bigl(\frac{2 }{(\lambda _{m}) ^2 L^4} \Bigr) ^{2M^*} t ^{8M^*}) e^{-t} \, dt 
\\
= 2^{2M^*}
\Bigl( 3^{2M^*} + \frac{2^{2M^*}\, (8M^*)! }{L^{8M^*} \, (\lambda _{m}) ^{4M^*}}\Bigr)  
\leq C_u ^{M^*}
\Bigl(1+
(\frac{1}{\lambda _m  (\gamma _{min} ^*)^2 (T')^2})^{4M^*} \, ({8M^*})! \Bigr)
.
\end{multline*}


Finally, we see that there exists some $C_u$ independent of $m$, $\gamma_{min}$, $\gamma _{min} ^*$, $N^*$ and $T$ such that
\begin{multline*}
\Vert \sigma _{m,N',T'} ^+ \Vert _{L^2(0,T)} ^2 
\leq C_u e^{-2\lambda _{m} T}
e^{C_u \frac{\sqrt{\lambda _m}}{\gamma _{min} ^*}} 
e^{C_u \frac{\lambda _{N^*}}{\gamma _{min} \sqrt{\lambda _m}}}
e^{\frac{C_u}{(\gamma _{min} ^* )^2 T'}}
(\frac{1}{(T')^{3/2}} + \frac{1}{(\gamma _{min} ^* )^2 (T')^2})
\\
e ^{C_u M^*}
\Bigl(1+
(\frac{1}{\lambda _{m}  (\gamma _{min} ^*)^2 (T')^2})^{4M^*} \, ({8M^*})! \Bigr),
\end{multline*}
which gives \eqref{*famillebi_qm-norme} and completes the proof of Lemma \ref{*lem-propf_m} and of Theorem \ref{thm-biortho1-gen}. \qed 


\section{Proof of Theorem \ref{thm-guichal-gen}}
\label{sec-guichal}

\subsection{A lower bound for any biorthogonal family} \hfill

Denote $E(\Lambda ,T)$ the smallest closed subspace of $L^2(0,T)$ containing the functions 
$$ \varepsilon _{\lambda _{n}}: s \in (0,T)\mapsto e^{-\lambda _{n}s}, \quad n\geq 1 .$$ 
It follows from \eqref{gap-max} that
$$ \sum _{n=1} ^\infty \frac{1}{\lambda _{n}} < \infty ,$$
and then it is well-known (\cite{Schwartz, Redheffer}) that $E(\Lambda ,T)$ is a proper subspace of $L^2(0,T)$. 
Moreover, given $m\geq 1$, denote $\Lambda _{ m} := (\lambda _{k})_{k\neq m}$, and $E(\Lambda _{m},T)$ the smallest closed subspace of $L^2(0,T)$ containing
the functions $\varepsilon _{\lambda _{k}}$, with $k\geq 1$ and $k\neq m$ (it does not include $\varepsilon _{\lambda _{m}}$). 
Then consider $p _{m}$ the orthogonal projection of  $\varepsilon _{\lambda _{m}}$ on $E(\Lambda _{m},T)$, and $d_{T,m}$ the distance between $\varepsilon _{\lambda _{m}}$ and $E(\Lambda _{ m},T)$: we have
\begin{equation}
\label{*def-distance}
d_{T,m}^2 = \inf _{p \in E(\Lambda _{m},T) } \Vert \varepsilon _{\lambda _{m}} - p\Vert _{L^2(0,T)} ^2 = \int _0 ^T ( e^{-\lambda _{m}s} - p _{m} (s)) ^2 \, ds .
\end{equation}
Then $\varepsilon _{\lambda _{m}} - p_{m}$ is orthogonal to $E(\Lambda _{ m},T)$ , which implies that
$$ \forall n\neq m, \quad \int _0 ^T (e^{-\lambda _{m}s} - p_{m} (s)) e^{-\lambda _{n}s} \, ds =0,$$
and 
\begin{multline*}
  \int _0 ^T (e^{-\lambda _{m}s} - p_{ m} (s)) e^{-\lambda _{m}s} \, ds \\
= \int _0 ^T (e^{-\lambda _{m}s} - p_{m} (s)) (e^{-\lambda _{m}s} - p_{m} (s)) \, ds = d_{T,m}^2 .
\end{multline*}
Hence consider 
\begin{equation}
\label{*def-optimalbi}
\sigma _{m} ^- (s) := \frac{e^{-\lambda _{m} s} - p _{m} (s)}{d_{T,m}^2 }:
\end{equation}
the sequence of functions $(\sigma _{m} ^-)_{m\geq 1}$ is a biorthogonal family for the set $(\varepsilon _{\lambda _{n}})_{n\geq 1}
= (e^{-\lambda _n t})_{n\geq 1}$ in $L^2(0,T)$.

Moreover it is optimal in the following sense: if $(\tilde \sigma _{m} ^-)_{m\geq 1}$ is another biorthogonal family for the set $(\varepsilon _{\lambda _{n}})_{n\geq 1}$ in $L^2(0,T)$, then for all $m\geq 1$, $\tilde \sigma _{m} ^- - \sigma _{m} ^-$ is orthogonal to all $\varepsilon _{\lambda _{n}}$, hence to $E(\Lambda,T)$, hence to $\sigma _{m} ^-$ since $\sigma _{m} ^- \in E(\Lambda ,T)$. Hence
$$ \Vert \tilde \sigma _{m} ^- \Vert _{L^2(0,T)} ^2 = \Vert \sigma _{m} ^- \Vert _{L^2(0,T)} ^2 + \Vert \tilde \sigma _{m} ^- - \sigma _{m} ^-\Vert _{L^2(0,T)} ^2 \geq \Vert \sigma _{m} ^- \Vert _{L^2(0,T)} ^2  .$$
Therefore 
\begin{equation}
\label{*jeudi8oct0}
\Vert \tilde \sigma _{m} ^- \Vert _{L^2(0,T)} \geq  \Vert \sigma _{m} ^- \Vert _{L^2(0,T)} =\frac{1}{d_{T,m} } .
\end{equation}
Hence $\frac{1}{d_{T,m} }$ is a lower bound of every biorthogonal sequence $(\tilde \sigma _{m} ^-)_{m\geq 1}$; and a bound from above for $d_{T,m}$ gives a bound from below for every biorthogonal sequence.

At last, we note that if the sequence of functions $(\tilde \sigma _{m} ^+)_{m\geq 1}$ is a biorthogonal family for the set $(e^{\lambda _n t})_{n\geq 1}$ in $L^2(0,T)$, then 
$$ \int _0 ^T \tilde \sigma _m ^+ (T-s) e^{\lambda _m T} e^{-\lambda _n s } \, ds = \delta _{mn} ,$$
hence $(\tilde \sigma _m ^+ (T-s) e^{\lambda _m T})_m$ is biorthogonal for the set $(e^{-\lambda _n t})_{n\geq 1}$ in $L^2(0,T)$. This implies that
\begin{equation}
\label{*dim26fev}
\Vert \tilde \sigma _{m} ^+ \Vert _{L^2(0,T)} \geq \frac{e^{-\lambda _m T}}{d_{T,m} } .
\end{equation}
Hence $\frac{e^{-\lambda _m T}}{d_{T,m} }$ is a lower bound of every biorthogonal sequence $(\tilde \sigma _{m} ^+)_{m\geq 1}$.
In the following (Lemma \ref{lem-dist-gapmax}), we provide a bound from above for $d_{T,m}$, that will give a bound from below for every biorthogonal sequence $(\tilde \sigma _{m} ^+)_{m\geq 1}$.


\subsection{A general result for sums of exponentials} \hfill

Clearly,
$$ d_{T,m} \leq \Vert e^{- \lambda _{m}s} - p(s) \Vert _{L^2 (0,T)} $$
for all $p \in E(\Lambda _{m},T)$.
The idea used in G\"uichal \cite{Guichal} is to chose a particular element $p \in E(\Lambda _{ m},T)$
in order to provide an upper bound of $d_{T,m}$.
The first thing to note is the following:
consider $M \geq m$ and 
$$ q(s) := \sum_{i=1} ^{M+1} A_i e^{- \lambda _{i}s} $$
with coefficients $A_1, \cdots , A_{M+1}$. Then $q \in E(\Lambda _{m},T)$ if and only if $A_m =0$, and when $A_m \neq 0$, then 
$$ \frac{1}{A_m} q(s)=  e^{- \lambda _{m} s}+ \sum_{i=1} ^{m-1} \frac{A_i}{A_m} e^{- \lambda _{i}s} + \sum_{i=m+1} ^{M+1} \frac{A_i}{A_m} e^{- \lambda _{i}s} ,$$
hence
\begin{equation}
\label{*distance-gene}
\Vert \frac{1}{A_m} q(s) \Vert _{L^2 (0,T)} \geq d_{T,m} .
\end{equation}
We will choose the coefficients $A_1, \cdots , A_{M+1}$ so that
$$ q(0)= q'(0)=q''(0)= \cdots = q^{(M-1)}(0)=0, \quad q^{(M)} (0)=1 .$$
The following lemma is essentially extracted from G\"uichal \cite{Guichal}:

\begin{Lemma}
\label{lem-guichal}
Consider $M\geq 0$, and $0 < \lambda _1 < \cdots < \lambda _{M+1}$.

a) There exist coefficients $A_1, \cdots , A_{M+1}$ so that the function $q$ defined by
$$ q(s) := \sum_{i=1} ^{M+1} A_i e^{- \lambda _i s} $$
satisfies
$$ \begin{cases} 
q(0)= 0 \\
q'(0)=0 \\
q''(0)= 0 \\
\vdots \\
q^{(M-1)}(0)=0 \\
q^{(M)} (0)=1 .
\end{cases}$$
The coefficients are given by the following formulas:
\begin{equation}
\label{eq-guichal-coeff}
\forall k\in \{1, \cdots, M+1 \}, \quad A_k = \frac{1}{\prod_{i=1, i\neq k}^{M+1} (\lambda _i - \lambda _k) } .
\end{equation}

b) With this choice of coefficients, we have
\begin{equation}
\label{eq-guichal-encadr}
\forall s >0, \quad 0 < q(s) \leq \frac{s^{M}}{M!} e^{-\lambda _1 s} .
\end{equation}
\end{Lemma}

The only difference with G\"uichal \cite{Guichal} is the estimate \eqref{eq-guichal-encadr} which is more precise than the one obtained in \cite{Guichal}, Lemma 4:
$$ \forall s >0, \quad 0 < q(s) < \frac{s^{M}}{M!}  .$$ 
In the following, we prove \eqref{eq-guichal-encadr}, and in a sake of completeness, we give the main arguments for part a)
of Lemma \ref{lem-guichal}.

\noindent {\it Proof of Lemma \ref{lem-guichal}.}

a) We write the linear system
$$\begin{cases}
0 = q(0) = \sum_{i=1} ^{M+1} A_i  \\
0 = q'(0) = \sum_{i=1} ^{M+1} - \lambda _i A_i  \\
0 = q''(0) = \sum_{i=1} ^{M+1} (- \lambda _i)^2 A_i  \\
\cdots \\
0 = q^{(M-1)} (0) = \sum_{i=1} ^{M+1} (- \lambda _i)^{(M-1)} A_i \\
1 = q^{(M)} (0) = \sum_{i=1} ^{M+1} (- \lambda _i)^M A_i .
\end{cases}$$
This can be written

\begin{equation}
\label{*systeme}
\left( \begin{array}{ccccc}
            1 & 1 & \cdots & \cdots & 1 \\
            -\lambda_1 & -\lambda _2 & \cdots & \cdots & - \lambda _{M+1} \\
            (-\lambda_1)^2 & (-\lambda _2)^2 & \cdots & \cdots & (- \lambda _{M+1})^2 \\
           \vdots & \vdots & \vdots & \vdots & \vdots \\
           (-\lambda_1)^M & (-\lambda _2)^M & \cdots & \cdots & (- \lambda _{M+1})^M 
\end{array} \right)
 \left( \begin{array}{c}
A_1 \\ A_2 \\ A_3 \\ \vdots \\ A_{M+1}  \end{array} \right)
=  \left( \begin{array}{c}
0 \\ 0 \\ 0 \\ \vdots \\ 1  \end{array} \right) .
\end{equation}
The $(M+1)\times (M+1)$ matrix ${\mathcal{A}}$ that appears in the left hand side of \eqref{*systeme} is invertible: indeed, its determinant is of Vandermonde type, and
$$ \det {\mathcal{A}} = \prod _{k<l}\Bigl( (-\lambda _l) - (- \lambda _k)\Bigr) =\prod _{k<l}\Bigl( \lambda _k -  \lambda _l \Bigr) \neq 0 .$$
Hence the system \eqref{*systeme} is invertible, and the Cramer's formula gives 
$$ A_k = \frac{\det {\mathcal{B}}}{\det {\mathcal{A}}},$$
where
${\mathcal{B}}$ is the $(M+1)\times (M+1)$ matrix obtained from ${\mathcal{A}}$ putting the right-hand side member of \eqref{*systeme} at the place of the $k^\text{th}$-column of ${\mathcal{A}}$. But then, we can develop $\det {\mathcal{B}}$ with respect to the $k^\text{th}$-column and we find again a Vandermonde determinant. Then using the formula of Vandermonde determinant, one gets \eqref{eq-guichal-coeff}.

b) We prove \eqref{eq-guichal-encadr} by induction. When $M=0$, \eqref{eq-guichal-encadr} is true. Assume that it is true for some $M$, and let us prove that it is true for $M+1$: take
$$ q(s) := \sum_{i=1} ^{M+2} A_i e^{- \lambda _i s} ,$$
where the coefficients $A_1, \cdots , A_{M+2}$, 	are chosen so that
$$ q(0)= q'(0)=q''(0)= \cdots = q^{(M)}(0)=0, \quad q^{(M+1)}(0)=1 .$$
Then the Taylor developments of $q$ and $q'$ say that
$$ q(s) = \frac{s^{M+1}}{(M+1)!} + O(s^{M+2}) \quad
\text{ and } \quad q'(s) = \frac{s^M}{M!} + O(s^{M+1})\quad \text{ as } s \to 0 .$$
Consider 
$$ \tilde q (s) = e^{-2\lambda _{M+2}s} \frac{d}{ds} ( q(s) e^{\lambda _{M+2}s}) .$$
Then
\begin{multline*}
\tilde q(s) = e^{-2\lambda _{M+2}s} \frac{d}{ds} \Bigl( \sum_{i=1} ^{M+2} A_i e^{(\lambda _{M+2}- \lambda _i) s} \Bigr) 
\\
= e^{-2\lambda _{M+2}s} \Bigl( \sum_{i=1} ^{M+2} A_i (\lambda _{M+2}- \lambda _i) e^{(\lambda _{M+2}- \lambda _i) s} \Bigr)
\\
= \sum_{i=1} ^{M+2} A_i (\lambda _{M+2}- \lambda _i) e^{-(\lambda _{M+2}+ \lambda _i) s}. 
\end{multline*}
But the last term in the series is clearly equal to $0$, hence $\tilde q$ is a sum of $M+1$ exponentials. Moreover,
\begin{multline*}
 \tilde q (s) = q'(s)  e^{-\lambda _{M+2}s} + \lambda _{M+2} q(s)  e^{-\lambda _{M+2}s}
\\
= (\frac{s^M}{M!} + O(s^{M+1}))   e^{-\lambda _{M+2}s} + \lambda _{M+2} (\frac{s^{M+1}}{(M+1)!} + O(s^{M+2})) e^{-\lambda _{M+2}s} 
\\
= \frac{s^M}{M!} + O(s^{M+1}) \quad \text{ as } s \to 0 .
\end{multline*}
Hence 
$$  \tilde q(0)= \tilde q'(0)= \tilde q''(0)= \cdots = \tilde q^{(M-1)}(0)=0, \quad \tilde q^{(M)}(0)=1 ,$$
and we can apply the induction assumption to $\tilde q$: then 
$$ 0 < \tilde q(s) < \frac{s^M}{M!} e^{-(\lambda _{M+2}+\lambda _1) s} .$$
We deduce first that $s \mapsto q(s) e^{\lambda _{M+2}s}$ is increasing. Since its value in $0$ is $0$, then $q$ is positive on $(0, +\infty)$. Next, we obtain that
\begin{multline*}
 \frac{d}{ds} ( q(s) e^{\lambda _{M+2}s}) \leq  \frac{s^M}{M!} e^{(2\lambda _{M+2}-\lambda _{M+2}- \lambda _1) s} 
\\
= \frac{s^M}{M!} e^{(\lambda _{M+2}- \lambda _1) s} 
\leq \frac{d}{ds}\Bigl( \frac{s^{M+1}}{(M+1)!} e^{(\lambda _{M+2}- \lambda _1) s} \Bigr) , 
\end{multline*}
hence by integration, 
$$q(s) e^{\lambda _{M+2}s} \leq \frac{s^{M+1}}{(M+1)!} e^{(\lambda _{M+2}- \lambda _1) s} ,$$
hence
$$q(s) \leq \frac{s^{M+1}}{(M+1)!} e^{- \lambda _1 s}, $$
which completes the induction argument and the proof of Lemma \ref{lem-guichal}. \qed 


\subsection{A precise estimate of the remaining part of the exponential function} \hfill

It turns out that we will need an estimate for the remaining part of the exponential function
$$ \sum _{n=N} ^\infty \frac{x^n}{n!} 
$$
in function of $x$ and $N$.
We prove the following general and precise result:

\begin{Lemma}
\label{l-reste-exp}
We have the following estimates:
\begin{equation}
\label{eq-estim-exp}
\forall N\geq 1, \forall x \geq 0, \quad 
\frac{1}{N!} \Bigl(\frac{x}{1+x} \Bigr) ^N e^x \leq \sum _{n=N} ^\infty \frac{x^n}{n!}
\leq C_1 N \Bigl(\frac{x}{1+x} \Bigr) ^N e^x ,
\end{equation}
where 
$$ C_1 = \max _{x\in \Bbb R_+ } \frac{(1-e^{-x}) (1+x)}{x} .$$

\end{Lemma}

\noindent {\it Proof of Lemma \ref{l-reste-exp}.} Denote
$$ f_N (x) := \sum _{n=N} ^\infty \frac{x^n}{n!}.$$
Let us prove by induction that
$$ \forall N\geq 0, \forall x \geq 0, \quad f_N (x) \geq  \frac{1}{N!} \Bigl(\frac{x}{1+x} \Bigr) ^N e^x .$$
First, of course $f_0(x)=e^x$, and then
$$ f_0 (x) \geq \frac{1}{0!} \Bigl(\frac{x}{1+x} \Bigr) ^0 e^x .$$
Next, assume that 
$$ \forall x \geq 0, \quad f_N (x) \geq \frac{1}{N!} \Bigl(\frac{x}{1+x} \Bigr) ^N e^x .$$
We note that $$ f_{N+1} '(x) = f_N (x) ,$$
and
$$ \frac{d}{dx} \Bigl(  \frac{1}{(N+1)!} (\frac{x}{1+x}) ^{N+1} e^x \Bigr) 
= \frac{1}{(N+1)!}  (\frac{x}{1+x}) ^{N} e^x \Bigl( \frac{N+1}{(1+x)^2} + \frac{x}{1+x} \Bigr)
.$$
The study of the variations of the function $x\mapsto \frac{N+1}{(1+x)^2} + \frac{x}{1+x} $
gives 
$$ \forall x\geq 0 , \quad 1-\frac{1}{4(N+1)} \leq \frac{N+1}{(1+x)^2} + \frac{x}{1+x} \leq N+1 ,$$
hence
$$ \frac{d}{dx} \Bigl(  \frac{1}{(N+1)!} (\frac{x}{1+x}) ^{N+1} e^x \Bigr) 
\leq \frac{N+1}{(N+1)!}  (\frac{x}{1+x}) ^{N} e^x = \frac{1}{N!}  (\frac{x}{1+x}) ^{N} e^x .$$
Then
$$ f_{N+1} '(x) = f_N (x) \geq \frac{d}{dx} \Bigl(  \frac{1}{(N+1)!} (\frac{x}{1+x}) ^{N+1} e^x \Bigr) ,$$
and since the values at $0$ are $0$, we obtain that
$$ \forall x \geq 0, \quad f_{N+1} (x) \geq \frac{1}{(N+1)!} (\frac{x}{1+x}) ^{N+1} e^x .$$
This proves the first part of \eqref{eq-estim-exp}.

For the second part (which is not necessary for us here), we note that 
$$ \forall x \geq 0, \quad f_1 (x) \leq C_1  \Bigl(\frac{x}{1+x} \Bigr) ^1 e^x .$$
Assume that 
$$ \forall x \geq 0, \quad f_N (x)  \leq C_1 N \Bigl(\frac{x}{1+x} \Bigr) ^N e^x .$$
Then
$$ \frac{d}{dx} \Bigl(  (\frac{x}{1+x}) ^{N+1} e^x \Bigr) 
=  (\frac{x}{1+x}) ^{N} e^x \Bigl( \frac{N+1}{(1+x)^2} + \frac{x}{1+x} \Bigr)
\geq (1-\frac{1}{4(N+1)}) (\frac{x}{1+x}) ^{N} e^x .$$
Hence
$$ f_{N+1} '(x) = f_N (x) \leq C_1 N \Bigl(\frac{x}{1+x} \Bigr) ^N e^x \leq \frac{C_1 N}{1-\frac{1}{4(N+1)}} \frac{d}{dx} \Bigl(  (\frac{x}{1+x}) ^{N+1} e^x \Bigr) .$$
To conclude, note that
$$ \forall N\geq 1, \quad \frac{N}{1-\frac{1}{4(N+1)}} \leq N+1:$$
indeed, 
$$ (N+1) (1-\frac{1}{4(N+1)}) = N+1-\frac{1}{4} = N + \frac{3}{4} \geq N .$$
Hence, we obtain that 
$$ \forall x \geq 0, \quad f_{N+1} (x)  \leq C_1 (N+1) \Bigl(\frac{x}{1+x} \Bigr) ^{N+1} e^x ,$$
which concludes the induction, and the proof of \eqref{eq-estim-exp}. \qed


\subsection{Consequence: a bound from above for the distance $d_{T,m}$} \hfill

As a consequence of the upper estimate \eqref{*distance-gene} for the distance and of Lemma \ref{lem-guichal}, we obtain the following inequality: for all $m \geq 1$, for all $M \geq m$, we have 

\begin{equation}
\label{*estim-dist}
d_{T,m}
 \leq \Bigl( \prod_{i=1, i\neq m}^{M+1} \vert \lambda _{i} - \lambda _{m} \vert \Bigr) \Bigl( \int _0 ^T \frac{s^{2M}}{M!^2} e^{-2 \lambda _{ 1} s}\, ds \Bigr) ^{1/2} 
.
\end{equation}
It remains to estimate the terms that appear in the right hand side.
This is the object of the next sections, and it is based on the gap conditions \eqref{gap-max} and \eqref{gap-max*}. 


\subsubsection{Estimate under the uniform gap condition \eqref{gap-max}} \hfill

We prove the following:

\begin{Lemma}
\label{lem-dist-gapmax-1gap}
Assume that $(\lambda _n)_n$ satisfies \eqref{gap-max}. 
Denote
$$ k_* := [\frac{2 \sqrt{\lambda _1}}{\gamma _{max}}] +m+2$$
and
$$ C(T, \gamma _{max}, \lambda _1,m) = 
\frac{6 \sqrt{1+2T\lambda _1}}{\pi ^2 \sqrt{2T}} \, 
\frac{(k_*-1)!}{(m+k_*+3)! \, (m-1)!}  \, 
\frac{( T \gamma _{max} ^2)^{k_*+2}}{(1+T \gamma _{max} ^2)^{m+k_*+3}} .$$
Then 
\begin{equation}
\label{min-distance-1gap}
\forall m \geq 1, \quad \frac{1}{ d_{T,m}} \geq C(T, \gamma _{max}, \lambda _1,m) \, e^{\frac{1}{T \gamma _{max} ^{2}}} .
\end{equation}
\end{Lemma}

\noindent {\it Proof of Lemma \ref{lem-dist-gapmax-1gap}.} 
Of course
$$  \int _0 ^T \frac{s^{2M}}{M!^2} e^{-2 \lambda _{1} s}\, ds \leq \frac{T ^{2M+1}}{M!^2 (2M+1)} ,$$
an, on the other hand,
$$  \int _0 ^T \frac{s^{2M}}{M!^2} e^{-2 \lambda _{1} s}\, ds
\leq \frac{T^{2M}}{M!^2} \int _0 ^T  e^{-2 \lambda _{ 1} s}\, ds
\leq  \frac{T^{2M}}{M!^2} \frac{1-e^{-2 \lambda _{ 1} T}}{2 \lambda _{ 1}} .$$
Hence 
$$ \int _0 ^T \frac{s^{2M}}{M!^2} e^{-2 \lambda _{1} s}\, ds
\leq \frac{T^{2M}}{M!^2} \inf \{\frac{T}{ 2M+1}, \frac{1-e^{-2 \lambda _{ 1} T}}{2 \lambda _{ 1}}  \} .$$
But it is easy to check that
$$ \forall a,b > 0, \quad \inf \{a,\frac{1}{b} \} \leq \frac{2a}{1+ab} .$$
Indeed, $\inf \{a,\frac{1}{b} \} = a$ if $ab \leq 1$, and in this case $1+ab\leq 2$, hence $a (1+ab) \leq 2a$. On the other hand, when $ab\geq 1$, $\inf \{a,\frac{1}{b} \} = \frac{1}{b}$, and $1+ab\leq 2ab$. We deduce that
\begin{equation}
\label{*estim-dist3}
\Bigl( \int _0 ^T \frac{s^{2M}}{M!^2} e^{-2 \lambda _{1} s}\, ds \Bigr)^{1/2}
\leq \frac{T^{M}}{M!} \frac{\sqrt{2T}}{\sqrt{2M+1+2T \lambda _{ 1}}} .
\end{equation}

Now it remains to estimate the product
$$ \prod_{i=1, i\neq m}^{M+1} \vert \lambda _{i} - \lambda _{m} \vert 
=  \Bigl( \prod_{i=1, i\neq m}^{M+1} \vert \sqrt{\lambda _{i}} - \sqrt{\lambda _{m}} \vert \Bigr)
\Bigl( \prod_{i=1, i\neq m}^{M+1} (\sqrt{\lambda _{i}} + \sqrt{\lambda _{m}} ) \Bigr).$$

We derive from \eqref{gap-max} first that

$$ \vert \sqrt{\lambda _{i}} - \sqrt{\lambda _{m}} \vert
\leq \gamma _{max} \vert i-m \vert ,$$
and next that
$$ \sqrt{\lambda _{i}} + \sqrt{\lambda _{m}}
\leq 2 \sqrt{\lambda _{1}} + \gamma _{max} (i+m) .$$
Hence 
\begin{itemize}
\item first
\begin{multline*}
\prod_{i=1, i\neq m}^{M+1} ( \sqrt{\lambda _{i}} + \sqrt{\lambda _{m}} )
\leq \prod_{i=1, i\neq m}^{M+1} (2 \sqrt{\lambda _{1}} + \gamma _{max} (i+m) )
\\
\leq \gamma _{max}^M \frac{(M+1+ [\frac{2 \sqrt{\lambda _1}}{\gamma _{max}}] +m+1)!}{([\frac{2 \sqrt{\lambda _1}}{\gamma _{max}}]+m+1)!}
= c^{(+)} \, \gamma _{max}^M \, (M+k_*)!
\end{multline*}
with
$$ c^{(+)} = \frac{1}{([\frac{2 \sqrt{\lambda _1}}{\gamma _{max}}]+m+1)!}
\quad \text{ and } \quad k_* := [\frac{2 \sqrt{\lambda _1}}{\gamma _{max}}] +m+2 ;$$
\item next
\begin{multline*}
\prod_{i=1, i\neq m}^{M+1} \vert \sqrt{\lambda _{i}} - \sqrt{\lambda _{m}} \vert
\leq \prod_{i=1, i\neq m}^{M+1} \gamma _{max} \vert i-m \vert
\\
= \gamma _{max} ^{M} (m-1)! (M-(m-1))!
= c^{(-)} \, \gamma _{max}^M \, (M-(m-1))! 
\end{multline*}
with 
$$ c^{(-)} = (m-1)! .$$
\end{itemize}
Combining this with \eqref{*estim-dist3}, we derive from \eqref{*estim-dist}
$$ d_{T,m} \leq c^{(+)} \, c^{(-)} \, \frac{\sqrt{2T} }{\sqrt{2M+1+2T \lambda _{ 1}}} 
(T \gamma _{max} ^2 )^M \frac{(M+k_*)! \, (M-m+1)! }{M!} .$$
Denote
$$ c_* := c^{(+)} \, c^{(-)} \, \frac{\sqrt{2T} }{\sqrt{1+2T \lambda _{ 1}}} . $$
Then, to conclude, we note that
\begin{multline*}
\frac{1}{d_{T,m}} = \frac{6}{\pi ^2} \sum _{M=m+1} ^\infty \frac{1}{(M-m)^2} \frac{1}{d_{T,m}}
\\
\geq \frac{6}{\pi ^2 c_* } \sum _{M=m+1} ^\infty \frac{1}{(M-m)^2} \frac{M!}{(M+k_*)! \, (M-m+1)! } (\frac{1}{T \gamma _{max} ^2}) ^M 
\\
\geq \frac{6}{\pi ^2 c_* } \sum _{M=m+1} ^\infty \frac{1}{(M+k_*+2)! } (\frac{1}{T \gamma _{max} ^2}) ^M 
\\
= \frac{6}{\pi ^2 c_* }( T \gamma _{max} ^2)^{k_*+2} \sum _{n=m+k_*+3} ^\infty \frac{1}{n! } (\frac{1}{T \gamma _{max} ^2}) ^n .
\end{multline*}
And using Lemma \ref{l-reste-exp}, we obtain that \eqref{min-distance-1gap}.
This gives the expected exponential behaviour in $1/(T \gamma _{max} ^2)$. In the following we take care of the asymptotic gap $\gamma _{max} ^*$.


\subsubsection{Estimate under the uniform gap condition \eqref{gap-max} and the asymptotic gap condition \eqref{gap-max*}} \hfill

Now, taking into account the "asymptotic gap" given by \eqref{gap-max*}, we will be able to improve the previous estimate, roughly speaking replacing $\gamma _{max} ^{2}$ by $(\gamma _{max} ^*) ^{2}$ in the exponential factor.

\begin{Lemma}
\label{lem-dist-gapmax} 
Assume that $(\lambda _n)_n$ satisfies \eqref{gap-max}-\eqref{gap-max*}. Then 

\begin{equation}
\label{def-coeffb}
\frac{1}{d_{T,m}} \geq  b^* (T,\gamma_{max},\gamma_{max}^*,N_*, \lambda _1,m) \, e^{\frac{1}{T (\gamma _{max} ^*)^2}}
\end{equation}
where $b^*$ is given by
\begin{itemize}
\item when $m\leq N_*$, we have
\begin{equation}
\label{min-distance}
b^*(T,\gamma_{max},\gamma_{max}^*,N_*, \lambda _1,m)
= C^* \frac{\sqrt{1+T\lambda _1}}{\sqrt{T}} \, \frac{(T \, (\gamma _{max} ^*) ^2)^{K_*+K' _*+2}}{(1+ (T \, (\gamma _{max} ^*) ^2))^{N_*+K_*+K' _*+3}} ,
\end{equation}
where 
$$ C^* = \frac{c_u (\gamma _{max} ^*)^{2(N_*-1)} }{C^{(+)} C^{(-)}} \frac{1}{(N_*+K_*+K'_*+3)!} ,$$
and $C^{(+)}$, $C^{(-)}$, $K_*$ and $K' _*$ are given respectively in \eqref{def-C+}, \eqref{def-C-}, \eqref{def-K*} and \eqref{def-K'*};

\item  when $m> N_*$, we have
\begin{equation}
\label{min-distance2}
b^* (T,\gamma_{max},\gamma_{max}^*,N_*, \lambda _1,m) 
= \tilde C ^* \frac{\sqrt{1+T\lambda _1}}{\sqrt{T}} \, \frac{( T (\gamma _{max} ^*)^2 )^{K_*+2}}{(1+T (\gamma _{max} ^*)^2)^{m+K_*+3}} ,
\end{equation}
where
$$ \tilde C ^* = \frac{c_u }{\tilde C^{(+)} \, \tilde C^{(-)}}  \frac{1}{(m+K_*+3)!} $$
where $\tilde C^{(+)}$, $\tilde C^{(-)}$ and $K_*$ are given respectively  in \eqref{def-tildeC-+}, \eqref{def-tildeC--} and \eqref{def-K*}.
\end{itemize}

\end{Lemma}

The starting point is of course \eqref{*estim-dist} and \eqref{*estim-dist3}. Concerning the estimate of the product, we proceed in the same way as previously, distinguishing several cases. 
%
%
We investigate what can be said when $m\leq N_* < M+1$: in this case, 
\begin{itemize}
\item first we see that
\begin{multline*}\forall i \geq N_*+1, \quad 
\sqrt{\lambda _{i}} + \sqrt{\lambda _{m}} = 
\sqrt{\lambda _{i}} - \sqrt{\lambda _{N_*}} + \sqrt{\lambda _{N_*}} + \sqrt{\lambda _{m}}
\\
\leq \gamma _{max} ^* (i-N_*) + 2 \sqrt{\lambda _1} + (m+N_*)\gamma _{max} ;
\end{multline*}
hence
\begin{multline}
\label{2-gap-+-inf}
\prod_{i=1, i\neq m}^{M+1} (\sqrt{\lambda _{i}} + \sqrt{\lambda _{m}} )
= \Bigl( \prod_{i=1, i\neq m}^{N_*} (\sqrt{\lambda _{i}} + \sqrt{\lambda _{m}} ) \Bigr)
\Bigl( \prod_{i=N_* +1}^{M+1} (\sqrt{\lambda _{i}} + \sqrt{\lambda _{m}} ) \Bigr)
\\
\leq \Bigl( \prod_{i=1, i\neq m}^{N_*} (2 \sqrt{\lambda _1} + \gamma_{max} (i+m)) \Bigr) 
\Bigl( \prod_{i=N_* +1}^{M+1} ( 2 \sqrt{\lambda _1} + (m+N_*)\gamma _{max} + \gamma _{max} ^* (i-N_*)) \Bigr) 
\\
\leq 
C^{(+)} (\gamma _{max} ^*) ^{M} \, (M+1-N_*+ [\frac{2\sqrt{\lambda _1}+(N_*+m)\gamma _{max}}{\gamma_{max}^*}]+1)! 
\\
= C^{(+)} (\gamma _{max} ^*) ^{M} \, (M+K_*)!
\end{multline}
with 
\begin{equation}
\label{def-C+}
C^{(+)} = (\frac{\gamma _{max}}{\gamma _{max} ^*}) ^{N_* -1} \frac{(N_*+m+ [\frac{2\sqrt{\lambda _1}}{\gamma_{max}}]+1)!}{(m+ [\frac{2\sqrt{\lambda _1}}{\gamma_{max}}]+1)! \, ([\frac{2\sqrt{\lambda _1}+(N_*+m)\gamma_{max}}{\gamma_{max}^*}]+1)! \, (2m+ [\frac{2\sqrt{\lambda _1}}{\gamma_{max}}]+1)} 
\end{equation}
and
\begin{equation}
\label{def-K*}
K_* := [\frac{2\sqrt{\lambda _1}+(N_*+m)\gamma _{max}}{\gamma_{max}^*}] - N_* + 2 ;
\end{equation}

\item next, similarly we have
\begin{multline}
\label{2-gap---inf}
\prod_{i=1, i\neq m}^{M+1} \vert \sqrt{\lambda _{i}} - \sqrt{\lambda _{m}} \vert
= \Bigl( \prod_{i=1, i\neq m}^{N_*} \vert \sqrt{\lambda _{i}} - \sqrt{\lambda _{m}} \vert \Bigr)
\Bigl( \prod_{i=N_* +1}^{M+1} \vert \sqrt{\lambda _{i}} - \sqrt{\lambda _{m}} \vert \Bigr)
\\
\leq \Bigl( \prod_{i=1, i\neq m}^{N_*} \gamma _{max} \vert i-m  \vert \Bigr)
\Bigl( \prod_{i=N_* +1}^{M+1} \gamma _{max} ^* (i-N_*) + \gamma _{max} (N_*-m) \Bigr)
\\
\leq C^{(-)} (\gamma _{max} ^*) ^{M} \, (M-N_*+2+ [\frac{\gamma_{max}}{\gamma_{max}^*}(N_*-m)])! 
\\
= C^{(-)} (\gamma _{max} ^*) ^{M} \, (M+K' _*)!
\end{multline}
with
\begin{equation}
\label{def-C-}
C^{(-)} =  (\frac{\gamma _{max}}{\gamma _{max} ^*}) ^{N_*-1} \, 
\frac{(m-1)! \, (N_*-m)! }{(1+ [\frac{\gamma_{max}}{\gamma_{max}^*}(N_*-m)])!} 
\end{equation}
and
\begin{equation}
\label{def-K'*}
K' _* :=  [\frac{\gamma_{max}}{\gamma_{max}^*}(N_*-m)] -N_*+2 ;
\end{equation}

\end{itemize}
We deduce from \eqref{*estim-dist}, \eqref{*estim-dist3}, \eqref{2-gap-+-inf} and \eqref{2-gap---sup} that

$$
d_{T,m} \leq C^{(+)} \, C^{(-)} \, \frac{\sqrt{2T}}{\sqrt{1+2T\lambda _1}} 
\, \frac{(M+K_*)! \, (M+K' _*)!}{M!} 
\, (T \, (\gamma _{max} ^*) ^2) ^M .$$
Denote
$$ C_* := C^{(+)} \, C^{(-)} \, \frac{\sqrt{2T}}{\sqrt{1+2T\lambda _1}} .$$
Hence
$$ d_{T,m} \leq C_* \, \frac{(M+K_*)! \, (M+K' _*)!}{M!} 
\, (T \, (\gamma _{max} ^*) ^2) ^M .$$
Then, as we did before, we have
\begin{multline*}
\frac{1}{d_{T,m}} = \frac{6}{\pi^2} \sum _{M=N_*+1} ^\infty \frac{1}{(M-N_*)^2} \frac{1}{d_{T,m}}
\\
\geq \frac{6}{\pi^2 \, C_*} \sum _{M=N_*+1} ^\infty \frac{1}{(M-N_*)^2} 
\frac{M!}{(M+K_*)! \, (M+K' _*)!} (\frac{1}{T \, (\gamma _{max} ^*) ^2} )^M .
\end{multline*}
Note that
\begin{multline*}
\frac{1}{(M-N_*)^2}  \frac{M!}{(M+K_*)! \, (M+K' _*)!}
= \frac{1}{(M-N_*)^2} \frac{1}{(M+K_*)! \, (M+1)\cdots (M+K' _*)}
\\
\geq  \frac{1}{(M+K_*+K' _* +2)!} .
\end{multline*}
Hence
\begin{multline}
\label{dist-avantder}
\frac{1}{d_{T,m}} 
\geq \frac{6}{\pi^2 \, C_*} \sum _{M=N_*+1} ^\infty \frac{1}{(M+K_*+K' _* +2)!}  (\frac{1}{T \, (\gamma _{max} ^*) ^2} )^M 
\\
= \frac{6}{\pi^2 \, C_*} (T \, (\gamma _{max} ^*) ^2 )^{K_*+K' _* +2} \sum _{n=N_*+K_*+K' _* + 3} ^\infty \frac{(\frac{1}{T \, (\gamma _{max} ^*) ^2} )^n}{n!}  .
\end{multline}
Applying Lemma \ref{l-reste-exp} to \eqref{dist-avantder}, we obtain
$$
\frac{1}{d_{T,m}} 
\geq \frac{6}{\pi^2 \, C_*} \frac{1}{(N_*+K_*+K' _*+3)!} \frac{X^{K_*+K' _*+2}}{(1+ X)^{N_*+K_*+ K' _* +3}}  e^{1/X} .
$$
with 
$$ X = T \, (\gamma _{max} ^*) ^2 .$$
This concludes the proof of Lemma \ref{lem-dist-gapmax} when $m\leq N*$. \qed

In the same way, if $m>N_*$, we have
\begin{itemize}

\item first
\begin{multline*}
\prod_{i=1, i\neq m}^{M+1} (\sqrt{\lambda _{i}} + \sqrt{\lambda _{m}} )
= \Bigl( \prod_{i=1}^{N_*} (\sqrt{\lambda _{i}} + \sqrt{\lambda _{m}} ) \Bigr)
\Bigl( \prod_{i=N_* +1, i\neq m}^{M+1} (\sqrt{\lambda _{i}} + \sqrt{\lambda _{m}} ) \Bigr)
\\
\leq \Bigl( \prod_{i=1}^{N_*} (2 \sqrt{\lambda _1} + \gamma_{max} (i+m)) \Bigr) 
\Bigl( \prod_{i=N_* +1, i\neq m}^{M+1} ( 2 \sqrt{\lambda _1} + (m+N_*)\gamma _{max} + \gamma _{max} ^* (i-N_*)) \Bigr) 
\\
\leq \Bigl( (\gamma _{max}) ^{N_*} \frac{(N_*+m+ [\frac{2\sqrt{\lambda _1}}{\gamma_{max}}]+1)!}{(m+ [\frac{2\sqrt{\lambda _1}}{\gamma_{max}}]+1)!}\Bigr)
\\
\Bigl( (\gamma _{max} ^*) ^{M-N_*} \frac{(M+1-N_*+ [\frac{2\sqrt{\lambda _1}+(N_*+m)\gamma _{max}}{\gamma_{max}^*}]+1)!}{([\frac{2\sqrt{\lambda _1}+(N_*+m)\gamma_{max}}{\gamma_{max}^*}]+1)! (m-N_*+ [\frac{2\sqrt{\lambda _1}+(N_*+m)\gamma _{max}}{\gamma_{max}^*}]+1)} \Bigr)
\\
= \tilde C ^{(+)} (\gamma _{max} ^*) ^{M} (M+K_*)!
\end{multline*}
with
\begin{multline}
\label{def-tildeC-+}
\tilde C ^{(+)} =  (\frac{\gamma _{max}}{\gamma _{max}^*} ) ^{N_*} \, 
\frac{(N_*+m+ [\frac{2\sqrt{\lambda _1}}{\gamma_{max}}]+1)!}{(m+ [\frac{2\sqrt{\lambda _1}}{\gamma_{max}}]+1)!}
\\
 \frac{1}{([\frac{2\sqrt{\lambda _1}+(N_*+m)\gamma_{max}}{\gamma_{max}^*}]+1)! (m-N_*+ [\frac{2\sqrt{\lambda _1}+(N_*+m)\gamma _{max}}{\gamma_{max}^*}]+1)}  ;
\end{multline}

\item next, 
\begin{multline}
\label{2-gap---sup}
\prod_{i=1, i\neq m}^{M+1} \vert \sqrt{\lambda _{i}} - \sqrt{\lambda _{m}} \vert
= \Bigl( \prod_{i=1}^{N_*} \vert \sqrt{\lambda _{i}} - \sqrt{\lambda _{m}} \vert \Bigr)
\Bigl( \prod_{i=N_* +1, i\neq m}^{M+1} \vert \sqrt{\lambda _{i}} - \sqrt{\lambda _{m}} \vert \Bigr)
\\
\leq \Bigl( \prod_{i=1}^{N_*} \gamma _{max} \vert i-m  \vert \Bigr)
\Bigl( \prod_{i=N_* +1, i\neq m}^{M+1} \gamma _{max} ^* \vert i-m \vert \Bigr)
\\
= (\gamma _{max}) ^{N_*} (\gamma _{max} ^*) ^{M-N_*}  \, (m-1)! \, (M+1-m)! 
\\
= \tilde C ^{(-)}\,  (\gamma _{max} ^*) ^{M} \, (M+1-m)!
\end{multline}
with
\begin{equation}
\label{def-tildeC--}
\tilde C ^{(-)} = (\frac{\gamma _{max}}{\gamma _{max} ^*} ) ^{N_*}  \, (m-1)! ;
\end{equation}

\item then we can conclude:

$$
d_{T,m} \leq \tilde C _* \frac{(M+K_*)! \, (M+1-m)! }{M!} \, T^M \, (\gamma _{max} ^*) ^{2M}
$$
with
\begin{equation}
\label{def-tildeC*}
\tilde C _* = \frac{\sqrt{2T}}{\sqrt{1+2T\lambda _1}}  \, \tilde C ^{(+)} \, \tilde C ^{(-)} ;
\end{equation}
then, in the same way, 
\begin{multline*}
\frac{1}{d_{T,m}} = \frac{6}{\pi ^2} \sum _{M=m+1} ^\infty \frac{1}{(M-m)^2} \frac{1}{d_{T,m}}
\\
\geq \frac{6}{\pi ^2 \tilde C _* } \sum _{M=m+1} ^\infty \frac{1}{(M-m)^2} \frac{M!}{(M+1-m)! \, (M+K_*)!} (\frac{1}{T (\gamma _{max} ^*)^2 })^M
\\
\geq \frac{6}{\pi ^2 \tilde C _* } \sum _{M=m+1} ^\infty \frac{1}{(M+K_*+2)!} (\frac{1}{T (\gamma _{max} ^*)^2 })^M
\\
= \frac{6}{\pi ^2 \tilde C _* } ( T (\gamma _{max} ^*)^2 )^{K_*+2} \sum _{n=m+K_*+3} ^\infty \frac{1}{n!} (\frac{1}{T (\gamma _{max} ^*)^2 })^n
\\
\geq \frac{6}{\pi ^2 \tilde C _* } ( T (\gamma _{max} ^*)^2 )^{K_*+2} \frac{1}{(m+K_*+3)!} \Bigl( \frac{1}{1+T (\gamma _{max} ^*)^2} \Bigr) ^{m+K_*+3} e^{1/(T (\gamma _{max} ^*)^2)} .
\end{multline*}
This concludes the proof of Lemma \ref{lem-dist-gapmax} when $m> N*$. \qed 

\end{itemize}



\begin{thebibliography}{99}
 



\bibitem{Assia1} F. Ammar Khodja, A. Benabdallah, M. Gonz\'alez-Burgos, L. de Teresa, {\it The
Kalman condition for the boundary controllability of coupled parabolic systems. Bounds on biorthogonal families to complex matrix exponentials}, J. Math. Pures Appl. 96 (2011), p. 555-590.

\bibitem{Karine-Grushin1} K. Beauchard, P. Cannarsa, R. Guglielmi,{\it Null controllability of Grushin-type operators in dimension two}, J. Eur. Math. Soc. (JEMS) 16 (2014), No. 1, 67-101.

\bibitem{Karine-Grushin2} K. Beauchard, L. Miller, M. Morancey, {\it 2D Grushin-type equations: minimal time and null controllable data}, J. Differential Equations. 259 (11), 2015.







\bibitem{sicon2008} 
P. Cannarsa, P. Martinez, J. Vancostenoble,
\textit{Carleman estimates for a class of degenerate parabolic operators},
SIAM J. Control Optim. 47, (2008), no. 1, 1--19.

\bibitem{memoire} P. Cannarsa, P. Martinez, J. Vancostenoble,
\textit{Global Carleman estimates for degenerate parabolic operators with applications}, Memoirs of the American Mathematical Society (2016), Vol. 239.

\bibitem{CMV-cost-weak} P. Cannarsa, P. Martinez, J. Vancostenoble,
{\it The cost of controlling weakly degenerate parabolic equations by boundary controls}, Math. Control Relat. Fields (2017), Vol 7, No 2, p. 171-211.


\bibitem{CMV-cost-loc} P. Cannarsa, P. Martinez, J. Vancostenoble,
{\it The cost of controlling strongly degenerate parabolic equations}, in preparation.




\bibitem{Coron} J.M. Coron, S. Guerrero, {\it Singular optimal control: A linear 1-D parabolic-hyperbolic example}, Asymp. Anal. 44, No 3-4, (2005), p. 237-257.


\bibitem{Jeremi-Sylvain} J. Dard\'e, S. Ervedoza, {\it On the reachable set for the one-dimensional heat equation} (2016), arXiv:1609.02692.


\bibitem{Ervedoza1} S. Ervedoza, E. Zuazua, {\it Sharp observability estimates for heat equations}, Arch. Ration. Mech. Anal. 202 (2011), No 3, 975-1017.





\bibitem{Everitt2} W.N. Everitt, {\it A catalogue of Sturm-Liouville differential equations}, Sturm-Liouville Theory, 271-331, Birkh\"auser, Basel (2005).

\bibitem{FR1} H. O. Fattorini, D. L. Russel, 
{\it Exact Controllability Theorems for Linear Parabolic Equations in One Space Dimension}, 
Arch. Rat. Mech. Anal. 4, 272-292 (1971).

\bibitem{FR2} H. O. Fattorini, D. L. Russel, 
{\it Uniform bounds on biorthogonal functions for real exponentials with an application to the control theory of parabolic equations}, Quart. Appl. Math. 32 (1974/75), 45-69.

\bibitem{FR3} H. O. Fattorini, {\it Boundary control of temperature distributions in a parallelepipedon}, SIAM J. Control, Vol 13, No 1 (1975).


\bibitem{FCEZ} E. Fernandez-Cara, E. Zuazua, 
{\it The cost of approximate controllability for heat equations: the linear case}, Adv. Differential equations 5 (2000), No 4-6, 465-514.


\bibitem{Fursikov}
A. V. Fursikov, O. Yu. Imanuvilov,
\textit{Controllability of evolution equations}. 
Lecture Notes Ser. 34, Seoul National University, Seoul, Korea, 1996.

\bibitem{Glass} O. Glass, {\it A complex-analytic approach to the problem of uniform controllability of  transport equation in the vanishing viscosity limit}, J. Funct. Anal. 258 (2010), No 3, p. 852-868.

\bibitem{Gueye} M. Gueye, {\it Exact boundary controllability of 1-D parabolic and hyperbolic degenerate equations}, SIAM J. Control Optim Vol 52 (2014), No 4, p. 2037-2054.

\bibitem{Guichal} E.N. G\"uichal, {\it A lower bound of the norm of the control operator for the heat equation}, Journal of Mathematical Analysis and Applications 110 (1985), p. 519-527.

\bibitem{Hansen} S. Hansen, {\it Bounds on functions biorthogonal to sets of complex exponentials; control of damped elastic systems}, Journal of Math. Anal. and Appl., 158 (1991), 487-508.

\bibitem{Haraux} A. Haraux, {\it S\'eries lacunaires et contr\^ole semi-interne des vibrations d’une plaque rectangulaire}, J. Math. Pures Appl. 68 (1989), 457-465.

\bibitem{Kamke} E. Kamke, Differentialgleichungen: L\"osungsmethoden und L\"osungen. Band 1: Gew\"ohnliche Differentialgleichungen. 3rd edition, Chelsea Publishing Company, New York, 1948. 


\bibitem{Kom-Lor} V. Komornik, P. Loreti, {\it Fourier Series in Control Theory}, Springer, Berlin, 2005.



\bibitem{Lions} J.-L. Lions, E. Zuazua, {\it On the cost of controlling unstable systems: the case of boundary controls}, J. Anal. Math. 73 (1997), 225-249.

\bibitem{Lissy1} P. Lissy, {\it On the cost of fast controls for some families of dispersive or parabolic equations in one space dimension}, SIAM J. Control Optim. 52 (2014), no. 4, 2651-2676.

\bibitem{Lissy2}  P. Lissy, {\it Explicit lower bounds for the cost of fast controls for some 1-D parabolic or dispersive equations, and a new lower bound concerning the uniform controllability of the 1-D transport-diffusion equation}, J. Differential Equations 259 (2015), no. 10, 5331-5352.


\bibitem{Lebedev}
N.N. Lebedev,
{\it Special Functions and their Applications},
Dover Publications, New York, 1972 

\bibitem{Rosier} 
P. Martin, L. Rosier, P. Rouchon, {\it Null controllability of one-dimensional parabolic equations using flatness}, Automatica J. IFAC 50 (2014) No 12, 3067-3076.

\bibitem{Rosier2} 
P. Martin, L. Rosier, P. Rouchon, {\it On the reachable states for the boundary control of the heat equation}, Applied Mathematics Research eXpress (2016).


\bibitem{Miller} L. Miller, {\it Geometric bounds on the growth rate of null controllability cost for the heat equation in small time}, J. Differential Equations 204 (2004), p. 202-226.





\bibitem{Redheffer} R.M. Redheffer, {\it Elementary remarks on completeness}, Duke Math. Journal 35 (1968), p. 103-116.

\bibitem{Schwartz} L. Schwartz, {\it \'Etude des sommes d'exponentielles}, deuxi\`eme \'edition. Paris, Hermann 1959.

\bibitem{Seidman} Th. Seidman, {\it Time invarinace of the reachable set for linear control problems}, J. Math. Annal. Appl. (1979), No 1, 17-20.

\bibitem{Seidman84} T. Seidman, {\it Two results on exact boundary control of parabolic equations}, Appl. Math. Optim. 11 (1984), p. 145-152.

\bibitem{Seid-Avdon} T.I. Seidman, S.A. Avdonin, S.A. Ivanov, {The "window problem" for series of complex exponentials}, J. Fourier Anal. Appl. 6 (2000),No 3, p. 233-254.

\bibitem{Tucsnak} G. Tenenbaum, M. Tucsnak, {\it New blow-up rates for fast controls of Schrodinger and heat equations}, J. Differential Equations 243 (2007), p. 70-100.

\bibitem{Watson} G. N. Watson, {\it A treatise on the theory of Bessel functions}, second edition,
Cambridge University Press, Cambridge, England, 1944. 

\bibitem{Young} R.M. Young, {\it An Introduction to Nonharmonic Fourier Series}, Academic Press, 1980.


\end{thebibliography}
\end{document}